\documentclass[preprint,12pt]{elsarticle_preprint}

\usepackage{amsmath}
\usepackage{amssymb}
\usepackage{amsopn}
\usepackage{amsfonts}
\usepackage{mathtools}
\usepackage[colorlinks = true,
            linkcolor = blue,
            urlcolor  = blue,
            citecolor = blue,
            anchorcolor = blue]{hyperref}

\usepackage[utf8]{inputenc}
\usepackage{subfig}

\usepackage{bm}
\usepackage{color}
\usepackage{comment}
\usepackage{graphicx}
\usepackage{subfig}
\usepackage{pgfplots}
\usepackage{todonotes}
\usepackage{multirow}
\usepackage{textcase}

\newcommand{\code}[1]{\texttt{#1}}

\renewcommand{\vec}[1]{\ensuremath{\bf #1}}

\newcommand{\trans}{{T}}

\newcommand{\vw}{{\vec{w}}}
\newcommand{\vc}{{\vec{c}}}
\newcommand{\vp}{{\vec{p}}}
\newcommand{\vb}{{\vec{b}}}
\newcommand{\vx}{{\vec{x}}}
\newcommand{\vu}{{\vec{u}}}
\newcommand{\vv}{{\vec{v}}}
\newcommand{\vf}{{\vec{f}}}

\newcommand{\vr}{{\vec{r}}}
\newcommand{\vy}{{\vec{y}}}
\newcommand{\vz}{{\vec{z}}}
\newcommand{\vzero}{{\vec{0}}}
\newcommand{\vmu}{{\vec{\mu}}}

\newcommand{\mI}{{\vec{I}}}
\newcommand{\mP}{{\vec{P}}}
\newcommand{\mE}{{\vec{E}}}
\newcommand{\mK}{{\vec{K}}}
\newcommand{\mA}{{\vec{A}}}

\newcommand{\mR}{{\vec{R}}}
\newcommand{\mM}{{\vec{M}}}
\newcommand{\mC}{{\vec{C}}}
\newcommand{\mW}{{\vec{W}}}
\newcommand{\mD}{{\vec{D}}}
\newcommand{\mQ}{{\vec{Q}}}
\newcommand{\mV}{{\vec{V}}}
\newcommand{\mPhi}{{\vec{\Phi}}}
\newcommand{\mPi}{{\vec{\Pi}}}
\newcommand{\mzero}{{\vec{0}}}
\newcommand{\precM}{\mM}

\newcommand{\deflspU}{\mathcal{U}}
\newcommand{\deflspV}{\mathcal{V}}

\DeclareMathOperator{\cov}{cov}
\providecommand{\cov}[1]{\cov\!#1}

\DeclareMathOperator{\var}{var}
\providecommand{\cov}[1]{\var\!#1}

\DeclareMathOperator{\expect}{\mathbb E}
\providecommand{\expect}[1]{\expect\!#1}

\DeclareMathOperator{\diam}{diam}
\providecommand{\diam}[1]{\diam\!#1}

\DeclareMathOperator{\Div}{Div}
\providecommand{\Div}[1]{\Div\!#1}

\newcommand{\KL}{\operatorname{\mathit{K\kern-.2em L}}}

\newcommand{\veks}[1]{\boldsymbol{#1}}      

\newcommand{\x}[0]{\boldsymbol{x}}

\newcommand{\U}[0]{\veks{u}}

\usepackage[normalem]{ulem}
\newcommand{\REMOVEDTEXTMODE}[1]{}
\newcommand{\REMOVED}[1]{
\ifmmode
\text{\REMOVEDTEXTMODE{$#1$}}
\else
\REMOVEDTEXTMODE{#1}
\fi}
\newcommand{\ADDEDTEXTMODE}[1]{{#1}}
\newcommand{\ADDED}[1]{
\ifmmode
\text{\ADDEDTEXTMODE{$#1$}}
\else
\ADDEDTEXTMODE{#1}
\fi}

\newcommand{\hl}[1]{{#1}}

\usepackage{pifont}

\usepackage{float}
\floatstyle{ruled}
\newfloat{algorithm}{tbp}{loa}
\providecommand{\algorithmname}{Algorithm}
\floatname{algorithm}{\protect\algorithmname}
\usepackage{algorithmic}

\begin{document}

\begin{frontmatter}

\title{Speeding up an unsteady flow simulation by adaptive BDDC and Krylov subspace recycling}

\author[1,2]{Martin Hanek}
\ead{martin.hanek@fs.cvut.cz}
\author[1]{Jan Papež}
\ead{papez@math.cas.cz}
\author[1]{Jakub Šístek\corref{cor1}}
\ead{sistek@math.cas.cz}
\cortext[cor1]{Corresponding author}

\address[1]{Institute of Mathematics of the Czech Academy of Sciences, {\v Z}itn{\' a} 25, 115 67 Prague, Czech Republic}
\address[2]{Faculty of Mechanical Engineering, Czech Technical University in Prague, Karlovo nám. 13, 121 35 Prague, Czech Republic}

\begin{abstract}
We deal with accelerating the solution of a sequence of large linear systems solved by preconditioned conjugate gradient method (PCG).
The sequence originates from time-stepping within a simulation of an unsteady incompressible flow.
We apply a pressure correction scheme and focus on the solution of the Poisson problem for the pressure corrector.
Its scalable solution presents the main computational challenge in many applications.
The right-hand side of the problem changes in each time step, while the system matrix is constant and symmetric positive definite.
The acceleration techniques are studied on a representative problem of flow around a unit sphere. 
Our baseline approach is based on a parallel solution of each problem in the sequence by nonoverlapping domain decomposition method.
The interface problem is solved by PCG with the three-level BDDC preconditioner.
As a preliminary step, an appropriate stopping criterion for the PCG iterations is chosen.
Next, two techniques for accelerating the solution are gradually added to the baseline approach.
Deflation is used within PCG with several approaches to Krylov subspace recycling.
Finally, we add the adaptive selection of the coarse space within the three-level BDDC method.
The paper is rich in experiments with careful measurements of computational times on a parallel supercomputer. 
The combination of the acceleration techniques eventually leads to saving more than 40~\% of the computational time.
\end{abstract}
	
	
\end{frontmatter}

\section{Introduction}

We study the problem of solving a sequence of linear systems with a constant matrix and variable right-hand sides.
There are many scenarios resulting in such sequences,
and we apply and study the methods on the Poisson problem of pressure (corrector) within an unsteady simulation of an incompressible flow.
As a representative model problem, we consider an unsteady flow around a unit sphere.

Unsteady incompressible flows of Newtonian fluids are modelled by the Navier-Stokes equations.
The system is time-dependent and nonlinear, and in its generality, it leads to solving a system of nonlinear equations repeatedly in each time step.
This approach can be very demanding with respect to computational resources,
and a number of widely used numerical schemes circumvent the need to solve these systems by splitting the whole problem into solving a sequence of 
simpler problems within each time step.
Let us mention the Pressure Implicit with Splitting of Operator (PISO) scheme~\cite{Issa-1986-SID},
implicit-explicit (IMEX) time integration~\cite{Karniadakis-1991-HOS},
or the pressure-correction methods (see, e.g., the review in \cite{guermondEtAl:06}),
to name a few.
In these approaches, a Poisson-type problem for pressure \hl{often} becomes the most time-consuming problem to solve,
and this trend worsens with the problem size, see, e.g., \cite{Sistek-2015-PIS}.
This is the reason why a considerable attention has been devoted in literature to accelerating the Poisson-type problems in this context
(e.g., \cite{Chentanez-2011-MFP,Costa-2018-FFT,Dick-2016-SFP,Golub-1998-FPS,mcadams-2010-PMP}).
To solve the problem in parallel, multigrid methods, preconditioned Krylov subspace methods, and combination thereof are often employed (see, e.g.,~\cite{Weller-1998-TAC}).

While using domain decomposition solvers for stand-alone Poisson problems is well established and studied in literature, 
their application in the context of sequences of Poisson problems arising in incompressible flow simulations is much less common.
We investigated the applicability of a nonoverlapping domain decomposition method for this task in~\cite{Hanek-2023-AMB},
considering the Balancing Domain Decomposition based on Constrains (BDDC) preconditioner introduced in~\cite{Dohrmann-2003-PSC}.
It provides the baseline approach also for the present paper.
More specifically, we use the three-level BDDC preconditioner~\cite{Tu-2007-TBT3D,Mandel-2008-MMB} for the preconditioned conjugate gradient (PCG) method.
Domain decomposition (DD) methods have a relatively expensive setup,
which includes factorization of the local subdomain matrices, and the matrix of the coarse problem.
It is then repeatedly used for each right-hand side in the time-stepping sequence.
The DD method reduces a global problem to the interface among subdomains,
and runs PCG on the Schur complement problem.
This leads to a significant reduction of the size of the vectors in the Krylov method.

We apply the finite element method (FEM) in connection with the incremental pressure-correction scheme to discretize the problem.
If the computational domain is fixed, the problems for the pressure corrector have a constant matrix and a new right-hand side vector in each time step.
Hence, the interface Schur complement problem and the BDDC preconditioner are set up only in the first time step and reused in the subsequent time steps.

In the rest of this section, we briefly review the techniques employed to accelerate the solution of the arising sequence of linear systems.
The gradual improvements eventually led us to saving more than 40 percent of the overall simulation time for realistic computations.
The techniques are then described in more detail in the subsequent sections.

\paragraph{Stopping criterion for Krylov iterations}

When solving a sequence of algebraic systems corresponding to time-dependent simulations with subsequent solutions close one to another, it seems appropriate to use the previously computed approximation as the starting guess for the current (new) system.
This however may not bring a significant improvement in the solution time unless an appropriate stopping criterion is used, as we will discuss and illustrate in numerical experiments.

In this study we terminate the iterations of PCG after a sufficient reduction of the relative residual.
However, we suggest to normalize the residual norm $\|\vr_k\|$ (where $\vr_{k} = \vb - \mA \vx_{k}$ is the residual associated with the approximation $\vx_{k}$ from the $k$-th iteration) by a quantity that is not related to the initial approximation.
For that purpose, we take the norm of the right-hand-side vector $\vb$ and terminate the iterations when
\begin{equation*}
\frac{\|\vr_{k}\|}{\|\vb\|} < 10^{-6}.
\end{equation*}
We compare this approach with normalizing by the norm of the initial residual $\vr_{0}$, that is stopping when
\begin{equation*}
\frac{\|\vr_{k}\|}{\|\vr_0\|} < 10^{-6},
\end{equation*}
which is inappropriate in the cases when a good initial approximation $\vx_{0}$ is available.

\paragraph{Krylov subspace recycling}
This technique, which involves reusing information from previous runs of a Krylov method, is frequently employed for nonsymmetric systems and restarted methods; see, e.g., \cite{Paretal06} or the survey~\cite{Soodhalter-2020-SRI}. Applying a subspace recycling to improve the convergence of the conjugate gradient (CG) method for solving a sequence of SPD linear systems has been considered, e.g., in \cite{Saad00,ErhGuy00}. 
Running the CG method in a subspace defined as the orthogonal complement of a stored basis is an established approach known as deflated CG~\cite{Dos88,Nic87,Saad00}.
Since CG typically reduces oscillating components of the error faster than smooth components,
making the Krylov subspace orthogonal to the eigenvectors corresponding to the smallest eigenvalues typically results in faster convergence.
As such, deflation can be seen as an alternative to preconditioning; see, e.g., \cite{Tanetal10} \hl{and \cite{KlaRhe12}}.
The deflation, on the other hand, requires a careful implementation as some of the operations in the approach are sensitive to accumulation of rounding errors and numerical loss of orthogonality.
Another drawback of these methods might be the relatively high memory cost of storing the deflation basis which limits its size.
Both of these issues are mitigated by the use of the BDDC preconditioner in our context.
This synergy was already realized for the earlier Finite Element Tearing and Interconnecting (FETI) method~\cite{Farhat-1991-MFE} applied to sequences of right-hand sides in~\cite{Farhat-1994-ESB}.
\hl{In~\cite{KlaRhe12}, the deflation is seen as an alternative way to implementing the coarse space in BDDC and FETI-DP methods.}

\paragraph{Adaptive selection of constraints in BDDC}
The BDDC preconditioner allows a flexible definition of its coarse space through defining the coarse degrees of freedom.
In particular, there are several variants of the adaptive BDDC (see the overview papers \cite{Klawonn-2016-ACS,Pechstein-2017-UFA}), 
and the common feature of the methods is solving a number of local eigenvalue problems,
followed by using selected eigenvectors for defining optimal coarse degrees of freedom.
The enriched coarse space typically improves the convergence of the iterative method
at the cost of more expensive preconditioner setup due to solving the eigenvalue problems.
While these approaches were developed mostly for problems with heterogeneous materials with largely varying coefficients,
they can be seen as a general approach to adjusting the strength of the DD preconditioners.
For this reason, they present an interesting option also for time-dependent problems studied in this paper,
in which any reduction of the number of iterations can lead to large savings of the computational time needed for the whole sequence.
In the present work, we apply the method described in~\cite{Mandel-2007-ASF,Mandel-2012-ABT}.
It is combined with the three-level BDDC preconditioner as a special case of adaptive multilevel BDDC~\cite{Sousedik-2013-AMB}.

\medskip
In our setting, we combine the adaptive BDDC preconditioner with deflated CG,
so that we can benefit from synergy of these approaches in several aspects.
Namely, the CG method runs on much shorter vectors (with the length given by the interface size rather than the global size) allowing storing more basis vectors of the deflation basis.
A nonstandard feature of combining BDDC with deflated CG is that the spectrum of the preconditioned operator is bounded from below by 1 with many eigenvalues clustered near this value.
This means that deflating eigenvectors corresponding to the smallest eigenvalues does not significantly improve the CG convergence.
Instead, we use the eigenvectors corresponding to the \emph{largest} eigenvalues of the preconditioned operator as the basis for deflation.
This can be understood as another level of preconditioning, further reducing the upper bound of the spectrum of the preconditioned operator; see, e.g., \cite{Tanetal10,KlaRhe12}.

In our case study (described in more details in Section~\ref{sec:numexp}),
the proper choice of stopping criterion can save 11~\% of the computational time spent on solving the system for pressure correction.
Considering the process with proper initial guess and stopping criterion as a baseline approach, the improvement of the subspace deflation is around 8~\%. Using adaptive coarse space in BDDC preconditioner saves nearly 19~\%. Finally, the combination of the subspace deflation and adaptive coarse space reduces the time by 23~\%.
The savings are even higher, around 43~\%, for a larger tested problem.
Nevertheless, the objective of this paper is not to demonstrate computational savings in a specific problem instance. It rather aims at presenting a set of techniques which, either individually or in combination, can be applied efficiently across a broad class of problems in incompressible flows.

\section{Model problem}

We consider a  domain~$\Omega \subset \mathbb{R}^3$ with its boundary $\partial\Omega$ consisting of three disjoint parts $\partial\Omega_S$, $\partial\Omega_\infty$, and $\partial\Omega_{O}$, 
$\partial\Omega = \partial\Omega_S \cup \partial\Omega_\infty \cup \partial\Omega_{O}$. 
Part $\partial\Omega_S$ is the interface between fluid and the rigid body, $\partial\Omega_\infty$ is the inflow free-stream boundary, and $\partial\Omega_{O}$ is the outflow boundary. 
The flow is governed by the Navier-Stokes equations of an incompressible viscous fluid, 
\begin{equation} 
	\label{eq:navierStokes}
	\begin{aligned}
  		\frac{\partial \U}{\partial t} + (\U \cdot \nabla ) \U - \nu \Delta \U + \nabla p &= \veks{0} \quad \mbox{in}\ \Omega , \\
 	 \nabla \cdot \U &= 0 \quad \mbox{in}\ \Omega ,
	\end{aligned}
\end{equation} 
where $\U(t,\x)$ is the velocity vector of the fluid, $t$ denotes time, $\nu$ is the kinematic viscosity of the fluid, and $p$ is the kinematic pressure. 
System~(\ref{eq:navierStokes}) is complemented by the initial and boundary conditions:
$\U (t=0, \x) = \veks{0}$ in $\Omega$,
$\U (t,   \x) = \U_{\infty}$ on $\partial\Omega_\infty$,
$\U (t,   \x) = \veks{0}$ on $\partial\Omega_S$, and
$ -\nu(\nabla \U){\veks{n}} + p{\veks{n}} = \veks{0}$ on $\partial\Omega_O$,
with $\veks{n}$ being the unit outer normal vector of $\partial\Omega$.

System~(\ref{eq:navierStokes}) can be efficiently solved by a pressure-correction method.
In particular, we use the incremental pressure-correction method in the \emph{rotational form} discussed by~\cite{guermondEtAl:06}.
Details of our implementation can be found in~\cite{Sistek-2015-PIS}.

In this approach, we first define the pressure increment (corrector) 
$\psi^{(n+1)} = p^{(n+1)} - p^{(n)}+ \nu \nabla \cdot \U^{(n+1)}$.
In order to compute the velocity and pressure fields $(\U^{(n+1)}, p^{(n+1)})$ at time $t^{(n+1)}$,
the following three subproblems are subsequently solved.
\begin{enumerate}
\item The velocity field $\U^{(n+1)}$ is obtained by solving the convection-diffusion problem for each component of velocity
\begin{equation}    
\label{eq:projection1} 
		\frac{1}{\Delta t} \U^{(n+1)} + (\U^{(n)} \cdot \nabla) \U^{(n+1)} - \nu \Delta 
   	 	\U^{(n+1)} =  \frac{1}{\Delta t} \U^{(n)}  - \nabla (p^{(n)} + \psi^{(n)})
    		\quad \mbox{in}\ \Omega
\end{equation}
for 
$\U^{(n+1)} = \U_{\infty}$ on $\partial\Omega_\infty$,
$\U^{(n+1)} = \veks{0}$ on $\partial\Omega_S$,
and $\nu(\nabla \U^{(n+1)}){\veks{n}} = p^{(n)}{\veks{n}}$ on $\partial\Omega_O$.

\item Next, the pressure corrector $\psi^{(n+1)}$ is obtained by solving the Poisson problem
\begin{equation}    
\label{eq:projection2}
    -\Delta \psi^{(n+1)} = - \frac{1}{\Delta t} \nabla\cdot \U^{(n+1)} \quad \mbox{in}\ \Omega
\end{equation}
for
$\frac{\partial \psi^{(n+1)}}{\partial \veks{n}} = 0$ on $\partial\Omega_\infty \cup \partial\Omega_S$
and $\psi^{(n+1)} = 0$ on $\partial\Omega_{O}$. 
\item
Finally, the  pressure field $p^{(n+1)}$ is updated with
\begin{equation}
  \label{eq:projection3}
  p^{(n+1)} = p^{(n)} + \psi^{(n+1)} - \nu \nabla \cdot \U^{(n+1)}.
\end{equation}
\end{enumerate}

Problems~(\ref{eq:projection1}), (\ref{eq:projection2}), and~(\ref{eq:projection3}) are discretized by the finite element method (FEM) using
Taylor-Hood $Q_2-Q_1$ hexahedral elements. 
They approximate the velocity and pressure fields by continuous piecewise tri-quadratic and tri-linear basis functions, respectively.
In the finite element mesh, there are $n_{\U}$ nodes with velocity unknowns and $n_{p}$ nodes with pressure unknowns, with the ratio $n_{\U} / n_{p}$ being approximately 8. 
\hl{Note that applying FEM to equation~(\ref{eq:projection3}) leads to solving an algebraic system corresponding to an $L_2$ projection.
The term $\nu \nabla \cdot \U^{(n+1)}$ on the right-hand side is discontinuous from one element to another, and hence it cannot be simply added to the pressure function.
}

For solving the algebraic problems arising from (\ref{eq:projection1}) and (\ref{eq:projection3}), we use the methods identified as optimal by \cite{Sistek-2015-PIS}.
In particular, the Generalized Minimal Residual method (GMRES) is used for (\ref{eq:projection1}), and the PCG method for~(\ref{eq:projection3}),
which corresponds to a linear solve related to an $L_2$-projection when FEM is applied.
Block Jacobi preconditioner using ILU(0) on the blocks corresponding to mesh partitions is used for both problems.

Problem~(\ref{eq:projection2}) translates to an algebraic system with a discrete Laplacian matrix of size $n_{p}\times n_{p}$ which is symmetric and positive definite for $\partial\Omega_{O} \ne \emptyset$,
i.e., a nonempty part with the `do-nothing' boundary condition,
\begin{equation}
   \label{eq:systems}
    \mK\, \vx^{(n+1)} = \vf^{(n+1)} \,. 
\end{equation}
This is a well-studied case from the point of view of domain decomposition methods, which are suitable solvers for this task.
The main focus of this study is a scalable solution of sequence~(\ref{eq:systems}) arising from the Poisson problem for pressure corrector (\ref{eq:projection2}).

\hl{In the case of $\partial\Omega_{O} = \emptyset$, problem~(\ref{eq:projection2}) corresponds to solving the so-called \emph{pure Neumann} problem,
and matrix $\mK$ is only semidefinite.
It has a one-dimensional nullspace of constant vectors corresponding to constant functions in the finite element space.
The deflated CG method can be efficiently employed to cope with this issue, by adding the constant vector into the deflation basis.}

\section{BDDC method with adaptive coarse space}
\label{sec:bddc}

As the baseline approach to solve the sequence of problems (\ref{eq:systems}),
we employ the three-level version~\cite{Tu-2007-TBT3D,Mandel-2008-MMB} of the BDDC method~\cite{Dohrmann-2003-PSC}
implemented in the
\code{BDDCML}\footnote{
\url{https://users.math.cas.cz/~sistek/software/bddcml.html}
}
library.
In particular, the approach is based on (i) reducing the global problem on the whole domain to the reduced problem defined at the interface between subdomains,
(ii) solving the reduced problem using PCG, while (iii) preconditioning the problem by the BDDC preconditioner.
For the sake of brevity, we drop the time index $(n+1)$ in the following discussion.
More precisely, vectors $\vx$ and $\vf$ without superscript are considered in the $(n+1)$-st time level,
while the superscript $(n)$ will be kept for their counterparts from the previous $n$-th time level.

\subsection{Iterative substructuring}
\label{sec:substructuring}

First, we consider the reduction of the global problem to the inter-sub\-do\-main interface. 
This procedure is rather standard and described, e.g., in monographs~\cite{Dolean-2015-IDD,Toselli-2005-DDM}.
To this end, the finite element mesh is divided into $N_S$ nonoverlapping subdomains $\Omega_i, i = 1,\dots,N_S$, with the partition respecting inter-element boundaries.
The subset of unknowns to which elements of subdomain $\Omega_i$ contribute is called local subdomain unknowns.

\hl{
In iterative substructuring,
local subdomain unknowns are further split to those belonging to just one subdomain, called \emph{interior} unknowns,
and the unknowns shared by several subdomains, which form the \emph{interface}~$\Gamma$.
Then we seek the solution of the global interface problem
\begin{equation}
	\label{eq:Sug}
	\mA \vx^{\Gamma} = \vb 
\end{equation}
using, for instance, the PCG method.
Here $\mA$ is the global Schur complement of the interior unknowns,
$\vx^{\Gamma}$ is the part of the solution vector corresponding to the interface,
and $\vb$ is the reduced right-hand side vector.
Note that the global Schur complement matrix $\mA$ is not explicitly constructed in iterative substructuring,
since only multiplications of vectors by $\mA$ are needed at each PCG iteration. 
}

\subsection{Multilevel BDDC preconditioner}
\label{sec:multilevelBDDC}

Multilevel BDDC preconditioner is used within PCG when solving the interface problem~(\ref{eq:Sug}).
More precisely, an action of the preconditioner $\mM_{BDDC}^{-1}$ produces a preconditioned residual $\vz^{\Gamma}$ from the residual in the $k$-th iteration 
$\vr^{\Gamma} = \mA \vx^{\Gamma}_{k} - \vb$
by implicitly solving the system $\mM_{BDDC} \, \vz^{\Gamma} = \vr^{\Gamma}$.

In the construction of the BDDC preconditioner, a set of coarse degrees of freedom required to be continuous among subdomains is selected.
If enough coarse degrees of freedom are defined for each subdomain,
the preconditioner corresponds to an invertible matrix.
In the baseline approach, we consider function values at selected interface nodes (corners) and arithmetic averages across subdomain faces and edges as the coarse degrees of freedom.
In adaptive BDDC, we further enrich this set by weighted averages over faces of subdomains derived from eigenvectors of generalized eigenvalue problems for each pair of subdomains sharing a face.
This approach is described in Section~\ref{sec:adaptiveBDDC}.

In the standard 2-level BDDC method, 
the coarse degrees of freedom define a global coarse problem with the unknowns $\vu_{C}$ and local subdomain problems with mutually independent degrees of freedom $\vu_{i}$.
\hl{How to obtain $\vu_{C}$ and $\vu_{i}$ will be discussed later.}

The BDDC method provides an approximate solution by combining the global coarse and local subdomain components as
\begin{equation} 
	\label{eq:averaging_subdomain_solves}
	\vz^{\Gamma} = \sum_{i=1}^{N_S} \left(\mR_{i}^{\Gamma}\right)^\trans \mD_{i} \mR_{Bi} \left( \vu_{i} + \mPhi_{i} \mR_{Ci} \vu_{C} \right),
\end{equation} 
\hl{where the restriction matrix $\mR_{i}^{\Gamma}$ selects the local interface unknowns from the global interface unknowns,
$\mR_{Bi}$ selects the local interface unknowns from those at the whole subdomain,}
the columns of $\mPhi_{i}$ contain the local coarse basis functions,
and $\mR_{Ci}$ is the restriction matrix of the global vector of coarse unknowns to those present at the $i$-th subdomain.
Matrix $\mD_i$ applies weights to satisfy the partition of unity.
In this work, it corresponds to a diagonal matrix with entries given either by the inverse cardinality of the set of subdomains sharing the interface unknown (denoted as \emph{card}) or derived from the diagonal entries of the local stiffness matrices (denoted as \emph{diag}).
These are rather standard choices in the DD literature; see, e.g.,~\cite{Toselli-2005-DDM}.

Vectors $\vu_{C}$ and $\vu_{i}$ in (\ref{eq:averaging_subdomain_solves}) are obtained in each iteration by solving
\begin{align} 
	\label{eq:coarse_problem}
	&\mK_{C} \vu_{C} = \sum_{i = 1}^{N_S} \mR_{Ci}^{\trans} \mPhi_{i}^{\trans} \mR_{Bi}^\trans \mD_i \mR_{i}^{\Gamma} \vr^{\Gamma}, \\
	\label{eq:localCorr}
	&\left[
	\begin{array}[c]{cc}
		\mK_{i} & \mC_{i}^{\trans} \\
		\mC_{i} & \mzero  
	\end{array}
	\right]
	\left[
	\begin{array}[c]{c}
		\vu_{i}\\
		\vmu_{i}
	\end{array}
	\right]
	=
	\left[
	\begin{array}[c]{c}
		\mR_{Bi}^\trans \mD_i \mR_{i}^{\Gamma} \vr^{\Gamma}\\
		\mzero
	\end{array}
	\right]\, , \, i=1,\dots,N_{S},
\end{align} 
where $\mK_{C}$ is the stiffness matrix of the global coarse problem,
$\mK_{i}$ is the local matrix assembled from elements in the $i$-th subdomain,
and $\mC_{i}$ is a constraint matrix enforcing zero values of the local coarse degrees of freedom in the second block row of~(\ref{eq:localCorr}).

When the number of subdomains reaches thousands, scalable solution of problem~(\ref{eq:coarse_problem}) by a direct solver becomes a challenge~\cite{Badia-2016-MBD,Sousedik-2013-AMB}.
A way to overcome this issue is to solve the coarse problem only approximately.
In multilevel BDDC~\cite{Tu-2007-TBT3D,Mandel-2008-MMB}, we apply the preconditioner to the coarse problem with subdomains playing the role of elements.
Details of our implementation within the \code{BDDCML} library can be found in~\cite{Sousedik-2013-AMB}.

As shown in~\cite{Mandel-2003-CBD}, the condition number of the Schur complement preconditioned by BDDC $\kappa(\mM_{BDDC}^{-1}\mA) = \lambda_{\max}/\lambda_{\min} \le C \left( 1 + \log ^2 (H/h) \right)$,
where $H$ is the characteristic subdomain size and $h$ is the characteristic element size.
In other words, the condition number is independent of the number of subdomains and grows only mildly with their size.
The analysis reveals that the smallest eigenvalue $\lambda_{\min} \ge 1$, and the bound for $\kappa(\mM_{BDDC}^{-1}\mA)$ is actually the upper bound on the largest eigenvalue $\lambda_{\max}$~\cite{Mandel-2005-ATP}.

\subsection{Adaptive selection of coarse degrees of freedom}
\label{sec:adaptiveBDDC}

The idea of adaptive BDDC is to enrich the coarse space with additional degrees of freedom.
These are chosen to reduce the condition number of the preconditioned operator $\kappa(\mM_{BDDC}^{-1}\mA)$ in a close-to-optimal way.
In the present work, we apply the adaptive multilevel BDDC method described in~\cite{Sousedik-2013-AMB}. 
Although adaptive BDDC was originally developed for problems requiring a large number of iterations, such as those with jumps in material parameters, the problem studied here represents another potential use case. In this case, strengthening the preconditioner, at the cost of additional setup computations, may reduce the overall number of iterations and, ultimately, the total computational time.

In adaptive BDDC, a number of generalized eigenvalue problems are solved, each corresponding to a pair of subdomains sharing a face.
Suppose that a face is shared by the $s$-th and the $t$-th subdomain.
The related generalized eigenvalue problem reads
\begin{equation}
\mPi\left( \mI-\mE_{st}\right) ^{\trans} \mA_{st} \left( \mI-\mE_{st}\right) \mPi \vw=\lambda\mPi \mA_{st}\mPi \vw,  \label{eq:eig-matrix}
\end{equation}
where $\mA_{st}$ is a block-diagonal matrix composed of the local Schur complements of the $s$-th and the $t$-th subdomain,
\begin{equation}
\mA_{st}=\left[ 
\begin{array}{cc}
\mA_{s} &  \\ 
& \mA_{t}%
\end{array}
\right],  
\label{eq:global-S}
\end{equation}
$\mPi$ is a projection matrix enforcing continuity of the coarse degrees of freedom initially defined on the subdomains,
such as arithmetic averages on edges and faces of the subdomains,
$\mI$ is the identity matrix,
and $\mE_{st}$ is the averaging matrix that makes the unknowns at the common interface continuous. 

Once the generalized eigenvalue problems are solved for each pair of subdomains sharing a face,
each eigenvector corresponding to an eigenvalue larger than a prescribed threshold $\tau$ is used to enrich the coarse space.
In particular, if $\lambda_{\ell} > \tau$,
its corresponding eigenvector $\vw_{\ell}$ is used to define a coarse degree of freedom $\vc^{st}_{\ell}$ as
\begin{equation}
\vc^{st}_{\ell}=\vw_{\ell}^{\trans}\mPi\left( \mI-\mE_{st}\right) ^{\trans} \mA_{st}\left( \mI-\mE_{st}\right) \mPi.
\label{eq:c-def}
\end{equation}
The part of the vector $\vc^{st}_{\ell}$ that corresponds to the unknowns within a face between the subdomains is added as a new row into the matrices $\mC_{s}$ and $\mC_{t}$ from the respective problems (\ref{eq:coarse_problem}).
In this work, arithmetic averages on edges and faces of the subdomains are used as the initial coarse degrees of freedom, and the eigenvectors corresponding to the largest eigenvalues of the generalized eigenvalue problems are used to enrich this set.
Details of the employed approach can be found in \cite{Sousedik-2013-AMB}.

The prescribed threshold $\tau > 1$ is an approximation of the target condition number of the preconditioned system $\kappa(\mM_{BDDC}^{-1}\mA)$.
With the minimal eigenvalue $\lambda_{\min} = 1$ given, the aim is to reduce the maximal eigenvalue $\lambda_{\max}$ according to $\tau$.
The smaller the threshold, the more coarse degrees of freedom are added.
This has the effect of reducing the number of iterations at the cost of increasing the size of the global coarse problem,
and hence the cost of each action of the BDDC preconditioner.
Consequently, there is a trade-off between pushing the number of iterations down and the cost of each of them in adaptive BDDC.
It should be also noted that although the solution of the local generalized eigenvalue problems is parallelized,
it still presents a significant overhead in the setup of the BDDC preconditioner that is realized before solving the problem in the first time step.

For implementation reasons, the number of eigenvectors computed on each face is limited to 10,
\hl{and the maximal number of iterations of the employed Locally Optimal Block Preconditioned Conjugate Gradient (LOBPCG) method~\cite{Knyazev-2001-TOP} is limited to 15.}
Due to these limits, the threshold $\tau$ is only a rather weak indicator of the final condition number of the preconditioned system,
typically underestimating it significantly.

\section{Deflation and Krylov subspace recycling}
\label{sec:deflation}

The principle of deflation for iterative solvers starts with decomposing the solution space $\mathbb{R}^n$ as $\mathbb{R}^n = \deflspU + \deflspV$, where $\deflspU$ is a \emph{deflation} subspace of a (relatively) small dimension, and $\deflspV$ is the orthogonal complement of $\deflspU$ with respect to a suitable inner product. Then the given problem is solved directly on $\deflspU$ and iteratively on $\deflspV$. If the subspace $\deflspU$ is chosen properly, one can get a significant speed-up of the iterative solution.

Deflation brings two challenges: i) how to construct $\deflspU$, and ii) how to implement the iterative method on $\deflspV$.
While the latter is a technical issue that has been resolved for many iterative (Krylov) methods and is independent of the choice of~$\deflspU$ (and $\deflspV$), a proper choice of the deflation subspace is problem-dependent and should be carefully addressed for each application.

In this section, we first recall the deflated preconditioned conjugate gradient method following the original literature \cite{Nic87,Dos88,Saad00}.
Then we present the choice of the deflation space that is, in our case, based on subspace recycling from~\cite{Saad00} with certain modifications.

\subsection{\hl{Deflated preconditioned conjugate gradient method}}

Given a symmetric positive definite matrix~$\mA \in \mathbb{R}^{n\times n}$, a right-hand side~$\vb \in \mathbb{R}^n$, and an initial approximation~$\vx_0$, define the initial residual $\vr_0 = \vb - \mA\vx_0$.
In the $j$-th iteration, CG (\cite{HesSti52}) generates an approximation~$\vx_j$ characterized by 
\[
    \vx_j \in \vx_0 + \mathcal{K}_j(\mA,\vr_0), \qquad \vr_j = \vb - \mA \vx_j \perp \mathcal{K}_j(\mA,\vr_0),
\]
where $\mathcal{K}_j(\mA,\vv) = \mbox{span}\{ \vv, \mA \vv, \ldots, \mA^{j-1}\vv\}$ is the $j$-th Krylov subspace.

Let columns of a matrix~$\mW$ give the basis of the subspace $\deflspU$. Then
\begin{equation}
    \label{eq:projectionQ}
    \mQ = \mI - \mW (\mW^T\mA\mW)^{-1} \mW^T \mA,
\end{equation}
is a projector on $\deflspV$, which is an orthogonal complement of $\deflspU$ with respect to the inner product induced by matrix~$\mA$.
For $\vv$ such that $\mW^T\vv = 0$ define
\[
    \hl{\mathcal{K}_{\deflspU,j}(\mA,\vv) \equiv \deflspU + \mathcal{K}_j(\mQ\mA\mQ,\vv).}
\]
From the properties of~$\mQ$, $\mathcal{K}_j(\mQ\mA\mQ,\vv) \cap \deflspU = \{0\}$.

Deflated Conjugate Gradient method is a modification of CG that generates approximations such that
\[
    \vx_j \in \vx_0 + \mathcal{K}_{\deflspU,j}(\mA,\vr_0), \qquad \vr_j = \vb - \mA \vx_j \perp \mathcal{K}_{\deflspU,j}(\mA,\vr_0),
\]
assuming that $\vx_0$ is such that $\mW^T \vr_0 = 0$.
A corresponding algorithm is given in \cite[Algorithm~3.5]{Saad00}; see also \cite[Sect.~3]{Dos88}.
For the ease of presentation, we will call~$\mW$ the deflation basis.

To speed up the convergence, a suitable preconditioning formally transforming the problem $\mA\vx=\vb$ into $\precM^{-1}\mA\vx = \precM^{-1}\vb$ is often considered.
Preconditioning can be also used in combination with the deflation. The corresponding algorithm is given, e.g., in \cite[Algorithm~3.6]{Saad00}, and we provide it below as Algorithm~\ref{alg:deflatedpcg} for completeness.

\begin{algorithm}[!ht]
\caption{Deflated PCG} \label{alg:deflatedpcg}

\begin{algorithmic}[1]
\STATE let $\mW = [\vw_1, \vw_2, \ldots, \vw_k]$ be a basis of~$\deflspU$
\STATE \textbf{input} $\mA$, $\vb$, preconditioner $\precM$, initial guess~$\vx_{-1}$
\STATE $\vx_0 = \vx_{-1} + \mW (\mW^T\mA\mW)^{-1} \mW^T (\vb - \mA \vx_{-1})$\hfill \textcolor{gray}{to assure that $\mW^T \vr_0 = 0$} \label{alg:firstprojection}
\STATE $\vr_{0}=\vb-\mA \vx_{0}$
\STATE $\vz_0 = \precM^{-1}\vr_0$
\STATE $\hat{\mu}_0 = (\mW^T\mA\mW)^{-1} \mW^T\mA \vz_0$ \label{alg:secondprojectiona}
\STATE $\vp_0 = \vz_0 - \mW\hat{\mu}_0$\hfill \textcolor{gray}{equivalently $\vp_0 = (\mI - \mW (\mW^T\mA\mW)^{-1} \mW^T \mA)\vz_0$} \label{alg:secondprojectionb}
\FOR{$j=1,2\dots$ until convergence}
    \STATE $\alpha_{j-1}= {\vr_{j-1}^T \vz_{j-1}}/{\vp_{j-1}^T \mA \vp_{j-1}}$
    
    \STATE $\vx_{j}=\vx_{j-1}+\alpha_{j-1}\vp_{j-1}$
    \STATE $\vr_{j}=\vr_{j-1}-\alpha_{j-1}\mA \vp_{j-1}$
    \STATE \textcolor{gray}{$\vr_{j} = \vr_{j}-\mW(\mW^T\mW)^{-1}\mW^T \vr_j$} \label{alg:reorthores} 
    \STATE $\vz_j = \precM^{-1} \vr_j$
    \STATE $\beta_{j-1}= {\vr_j^T \vz_j} / {\vr_{j-1}^T \vz_{j-1}}$
    \STATE $\hat{\mu}_j = (\mW^T\mA\mW)^{-1} \mW^T\mA \vz_j$ \label{alg:thirdprojectiona}
    \STATE $\vp_{j}= \vz_{j} + \beta_{j-1}\vp_{j-1} - \mW \hat{\mu}_j$ \\~\hfill \textcolor{gray}{$\vp_j = \beta_{j-1}\vp_{j-1} + (\mI - \mW (\mW^T\mA\mW)^{-1} \mW^T \mA)\vz_j$} \label{alg:thirdprojectionb}
    
    \smallskip
\ENDFOR
\end{algorithmic}

\end{algorithm}

In comparison with the standard PCG, deflated PCG additionally requires a computation of initial vector~$\vx_0$ from $\vx_{-1}$ (line~\ref{alg:firstprojection} of Algorithm~\ref{alg:deflatedpcg}) and a projection of $\vz_0$ onto~$\deflspV$, $\vp_0 = \mQ \vz_0 = (\mI - \mW (\mW^T\mA\mW)^{-1} \mW^T \mA)\vz_0$ (lines~\ref{alg:secondprojectiona} and \ref{alg:secondprojectionb} of Algorithm~\ref{alg:deflatedpcg}) in the initial setup.
These operations require solution of a system with the matrix $\mW^T\mA\mW$,
which has the size equal to the dimension of the deflation space.
Then another projection $\mQ\vz_j = (\mI - \mW (\mW^T\mA\mW)^{-1} \mW^T \mA)\vz_j$ must be additionally computed in each iteration (lines~\ref{alg:thirdprojectiona} and \ref{alg:thirdprojectionb} of Algorithm~\ref{alg:deflatedpcg}).
These operations and the need of storing the basis~$\mW$ increase the time and memory requirements of deflated PCG.
Hence, similarly to adaptive BDDC, the deflation leads to a trade-off between reducing the number of iterations and increasing their cost.

Mathematically, it holds that $\mW^T \vr_j = \mW^T \vp_j = 0$. In finite-precision computations, however, this may not hold due to the loss of orthogonality.
In such cases, the reorthogonalization of residuals at line~\ref{alg:reorthores} of Algorithm~\ref{alg:deflatedpcg} may be necessary; see also~\cite[Eq.~(7.1)]{Saad00}.

\hl{The deflation can be elegantly employed also in the case when matrix $\mA$ is only symmetric positive semidefinite with a known nullspace.
This is the case when problem (\ref{eq:projection2}) becomes a pure Neumann problem for $\partial\Omega_{O} = \emptyset$.
In this case, the nullspace of $\mA$ is spanned by constant vectors, and a basis vector of ones can be added as one column of $\mW$.
Note that the nullspace of Schur complement $\mA$ is a restriction of the nullspace of $\mK$ in (\ref{eq:systems}) from all unknowns to the interface $\Gamma$.
}

\subsection{Subspace recycling}
\label{sec:recycling}

To construct a deflation space~$\deflspU$, one can use some a~priori information on the eigenvector spaces corresponding to a problematic part of the spectrum, as was considered, e.g., in \cite{Dos88}.
However, we have no such information in hand for the considered application.

When solving a sequence of systems $\mA\vx^{(i)} = \vb^{(i)}$, the idea might be to reuse, or \emph{recycle}, the Krylov subspace built for the previous system(s) to construct the deflation space for the current system,
as in~\cite{ErhGuy00}.
We denote the vectors computed in (deflated) PCG applied to $\mA\vx^{(i)} = \vb^{(i)}$ by the superscript, for example $\vp^{(i)}_j$ denotes the $j$-th search vector for the $i$-th system.
The maximal dimension of the deflation space is a parameter to be chosen by the user, we denote it by~$R$.

In the experiments, we first test two simple constructions of the deflation basis:

\paragraph{B1 (first $R$ search vectors)} 
The deflation basis is set as the first $R$ search vectors and then remains unchanged for the rest of the computation, $\mW^{(i)} = \mW = [\vp^{(1)}_1, \vp^{(1)}_2, \ldots, \vp^{(1)}_R]$.
If~$R$ is larger than the number~$\ell$ of PCG iterations for solving the first system, the basis~$\mW$ also contains the search vectors from solving the second system (and possibly more systems until it has~$R$ columns), $\mW =[\vp^{(1)}_1, \ldots, \vp^{(1)}_\ell, \vp^{(2)}_1, \ldots]$.

\paragraph{B2 (last $R$ search vectors)}
The basis is first constructed as in the previous case, and then it uses a `sliding window' to contain the last $R$ search vectors from the previous systems.
This choice is motivated by the fact that the search vectors from the preceding linear systems may be more relevant to the subsequent system to be solved.
In more detail, after solving the $(i-1)$-th system and saving $\ell$~search vectors as columns of~$\mP^{(i-1)}$, the deflation basis is updated as
    \[
        \mW^{(i)} = \left[ \vw^{(i-1)}_{\ell+1}, \vw^{(i-1)}_{\ell+2}, \ldots, \vw^{(i-1)}_{R}, \vp^{(i-1)}_1, \vp^{(i-1)}_2, \ldots, \vp^{(i-1)}_\ell\right].
    \]
In other words, $\ell$ `oldest' vectors in $\mW^{(i-1)}$ are replaced by the most recent search vectors while keeping the number of columns of~$\mW$ equal to~$R$.

The search vectors~$\vp_j$ are in PCG, as well as in deflated PCG (\cite[Prop.~3.3]{Saad00}), $\mA$-orthogonal, $\vp_j^T\mA\vp_l = 0$ for $j \neq l$.
Therefore, if the basis~$\mW$ contains the search vectors only, the matrix $\mW^T\mA\mW$ is diagonal and the projection~$\mQ$ in (\ref{eq:projectionQ}) only consists of matrix multiplications.
In finite-precision computations, the $\mA$-orthogonality between the search vectors is often lost during the course of iterations.
However, in our application with the BDDC preconditioner, only a decent number of PCG iteration was needed to reach the tolerance, and the loss of orthogonality was modest in our experiments.

The choice of the deflation basis according to \emph{B1} and \emph{B2} allows for a cheaper (faster) iterations, but the overall speed-up can be higher for the more elaborated construction of the deflation space from~\cite{Saad00}.
In particular, the construction is based on storing (some of) the search direction vectors and determining the new basis~$\mW^{(i)}$ as the solution of an eigenvalue problem.
Note that such a basis is not necessarily $\mA$-orthogonal, and consequently, matrix $\mW^T\mA\mW$ is no longer diagonal.
This recycling procedure is as follows:

\begin{enumerate}
    \item Solve the first system $\mA\vx^{(1)} = \vb^{(1)}$.
          Save (some of) the search vectors~$\vp^{(1)}_j$ computed within the PCG iterations into the matrix $\mP^{(1)}$ and the associated vectors $\mA\vp^{(1)}_j$ into~$\mA\mP^{(1)}$.
          The maximum number of vectors to save is, in general, one of the parameters to be determined by the user.
          We use all search vectors from each PCG run in this paper.
          Set $\mW^{(1)} = \emptyset$.
    \item For the other systems in the sequence, $i=2,3, \ldots, n_i$:
    \begin{enumerate}
        \item Use vectors in $\mW^{(i-1)}$, $\mP^{(i-1)}$, and $\mA\mP^{(i-1)}$ to generate a generalized eigenvalue problem specified below.
        \item Solve the eigenvalue problem and determine matrix~$\mW^{(i)}$ using a part of the computed eigenvectors. These vectors ideally approximate some eigenvectors of~$\precM^{-1}\mA$.
        \item Solve $\mA\vx^{(i)} = \vb^{(i)}$ by deflated PCG with an initial guess~$\vx^{(i)}_{-1}$, preconditioner~$\precM^{-1}$, and the deflation space basis~$\mW^{(i)}$. As for the first system, construct successively the matrices $\mP^{(i)}$ and~$\mA\mP^{(i)}$.
    \end{enumerate}
\end{enumerate}

It remains to detail how the generalized eigenvalue problem is constructed.
The information in $\mW^{(i-1)}$ and $\mP^{(i-1)}$ is combined to improve the approximation to the eigenvectors of the preconditioned matrix $\precM^{-1}\mA$.
Denote by $\mV$ the concatenated matrix $\mV = [\mW^{(i-1)}, \mP^{(i-1)}]$, then the \emph{Ritz approximation} to the eigenpairs of~$\precM^{-1}\mA$ is given by solving
\begin{equation}
\label{eq:Ritz}
    \mV^T\mA\mV \vy = \theta \mV^T\precM\mV\vy
\end{equation}
and setting $\vw=\mV\vy$ as an approximation to an eigenvector of~$\precM^{-1}\mA$. This standard approximation, however, requires to apply the preconditioner~$\precM$ to~$\mV$, which may not be possible when only the operation $\precM^{-1}\vv$ is available, as in our case of the BDDC preconditioner.
As in~\cite[Sect.~5.1]{Saad00}, we therefore consider the \emph{harmonic Ritz approximation}
\begin{equation}
\label{eq:harmonicRitz}
    \mV^T \mA\precM^{-1}\mA\mV\vy = \theta \mV^T \mA \mV \vy
\end{equation}
that involves the operation with the inverse of the preconditioner.
Moreover, we employ a strategy from \cite[Sects.~5.1 and 5.2]{Saad00} for constructing matrices $\mV^T \mA\precM^{-1}\mA\mV$ and $\mV^T \mA \mV$ from \eqref{eq:harmonicRitz} for $\mV = [\mW^{(i-1)}, \mP^{(i-1)}]$ at a low cost.
Notice that while problems \eqref{eq:Ritz} and \eqref{eq:harmonicRitz} are different,
they both provide approximations to the eigenvalues and eigenvectors of the preconditioned system matrix~$\precM^{-1}\mA$.

After solving one or several systems, the maximal size~$R$ of the deflation basis~$\mW^{(i-1)}$ is reached, and only a subset of $R$ harmonic Ritz approximations (eigenpairs satisfying~\eqref{eq:harmonicRitz} with $\vw=\mV\vy$) are taken to form~$\mW^{(i)}$.
We consider two variants:

\paragraph{B3 (eliminating $R$ smallest Ritz values)} For $\mV = [\mW^{(i-1)}, \mP^{(i-1)}]$ and $\vy_l$ the generalized eigenvector from~\eqref{eq:harmonicRitz} corresponding to the $l$-th \emph{smallest} Ritz value~$\theta_l$, set the deflation vector $\vw_l = \mV\vy_l$, $l = 1, 2, \ldots, R$. This choice from~\cite{Saad00} is motivated by the fact that the smallest eigenvalues of~$\precM^{-1}\mA$ typically harm the convergence of (deflated) PCG the most.

\paragraph{B4 (eliminating $R$ largest Ritz values)} Analogous to \emph{B3}, but now the eigenvectors corresponding to the \emph{largest} Ritz values~$\theta_l$ are taken to construct~$\vw_l$.
This is a key modification for our application.
Since the BDDC preconditioner shifts the lower part of the spectrum to one, creating a cluster of eigenvalues there,
the deflation space would need to be very large to guarantee a faster convergence for~\emph{B3}.
Despite the fact that also the upper part of the spectrum is close to one (the largest eigenvalues were $O(1)$ in our experiments) and the largest eigenvalues and corresponding eigenspaces are typically implicitly well-approximated by Ritz values and vectors in PCG,
this choice for~$\mW^{(i)}$ gives an interesting speed-up. 

Finally, note that after solving some number of systems $\mA\vx^{(i)}=\vb^{(i)}$, the eigenpairs of~$\precM^{-1}\mA$ can be well approximated, and no further improvement is gained by solving~\eqref{eq:harmonicRitz}.
Therefore, we set a heuristic criterion
\begin{equation}
\label{eq:harmritzstopcrit}
    \frac{\| \boldsymbol{\theta}^{(i)} - \boldsymbol{\theta}^{(i-1)} \|_2}{\| \boldsymbol{\theta}^{(i)} \|_2 } \leq 10^{-5},
\end{equation}
where $\boldsymbol{\theta}^{(i)}$ denotes $R$~harmonic Ritz values given by \eqref{eq:harmonicRitz} for the $i$-th system and selected according to \textit{B3} or \textit{B4} described above.
When criterion (\ref{eq:harmritzstopcrit}) is satisfied, say after solving the $m$-th system, the deflation space approximates the space spanned by the $R$ eigenvectors corresponding to the largest (or smallest) eigenvalues sufficiently well, and the basis is fixed for all the subsequent systems, \mbox{$\mW^{(n)}=\mW^{(m)}$}, $n \geq m$.
Criterion~\eqref{eq:harmritzstopcrit} can be replaced by a more elaborated stopping criteria; see, e.g., \cite{BenBra94}. However, \eqref{eq:harmritzstopcrit} performs satisfactorily in our computations.

\hl{
An interesting approach for selecting the deflation basis with the aid of proper orthogonal decomposition (POD) was presented in~\cite{DiazCortes-2021-ASL} in the context of evolving system matrix.
The relation between this POD-based deflation and harmonic Ritz vectors, which are studied in the present paper, remains to be investigated.
}

\section{Stopping criteria based on (algebraic) residual}
\label{sec:residuals}

A common (often the default) criterion to terminate an iterative solver is a prescribed reduction of the relative residual as
\begin{equation}
\label{eq:stopping_residual}
\frac{\|\vr_{k}\|}{\|\vr_{0}\|} < tol,
\end{equation}
where $\vr_{k}$ is the residual in the $k$-th iteration, $\vr_{k} = \vb - \mA \vx_{k}$, and $\vr_{0} = \vb - \mA \vx_{0}$ is the residual associated with the initial guess~$\vx_{0}$.
Clearly, this criterion involves only (cheaply) computable quantities. On the other hand, one should be aware of its significant limitations; \eqref{eq:stopping_residual} does not guarantee a small norm of the error, the norm $\| \vr_{k} \|$ depends on the discretization basis, and it is not clear how to choose the tolerance $tol$.

There is another good reason why~\eqref{eq:stopping_residual} may be a bad choice when solving sequences of linear systems. 
An accurate approximation \hl{$\vx_{0} \approx \vx$ typic}ally makes the norm of the initial residual $\vr_0$ small. 
\hl{Such initial approximation might be given as the (final) approximation computed for the previous system or using more elaborated strategies, e.g., based on reduced order models and proper orthogonal decomposition as in~\cite{MarJan06}.}
The requirement of further improving the norm of the residual $\vr_k$ by several orders of magnitude then may lead to useless iterations and an effect of \emph{over-solving}; after certain level, the error in computing $\psi^{(n+1)}$ in \eqref{eq:projection2} is dominated by the FEM discretization error and reducing the algebraic error further brings no improvement. 
Consequently, it is important to normalize $\|\vr_k\|$ by a quantity that is not related to the initial approximation.
In this study, we take the norm of the right-hand side vector $\vb$ for this purpose, and we terminate the iterations when
\begin{equation}
\label{eq:stopping_rhs}
\frac{\|\vr_{k}\|}{\|\vb\|} < 10^{-6}.
\end{equation}

Finally,
note that \eqref{eq:stopping_rhs} is used as a stopping criterion for PCG for example in MATLAB~R2019b, SciPy (all versions), or PETSc~3.21.
However, \eqref{eq:stopping_residual} was the default option in PETSc until version~3.18, and it might still be used in some other software packages.

\section{Results}
\label{sec:numexp}

\begin{figure}[htbp]
\centering
    \includegraphics[width=0.9\textwidth]{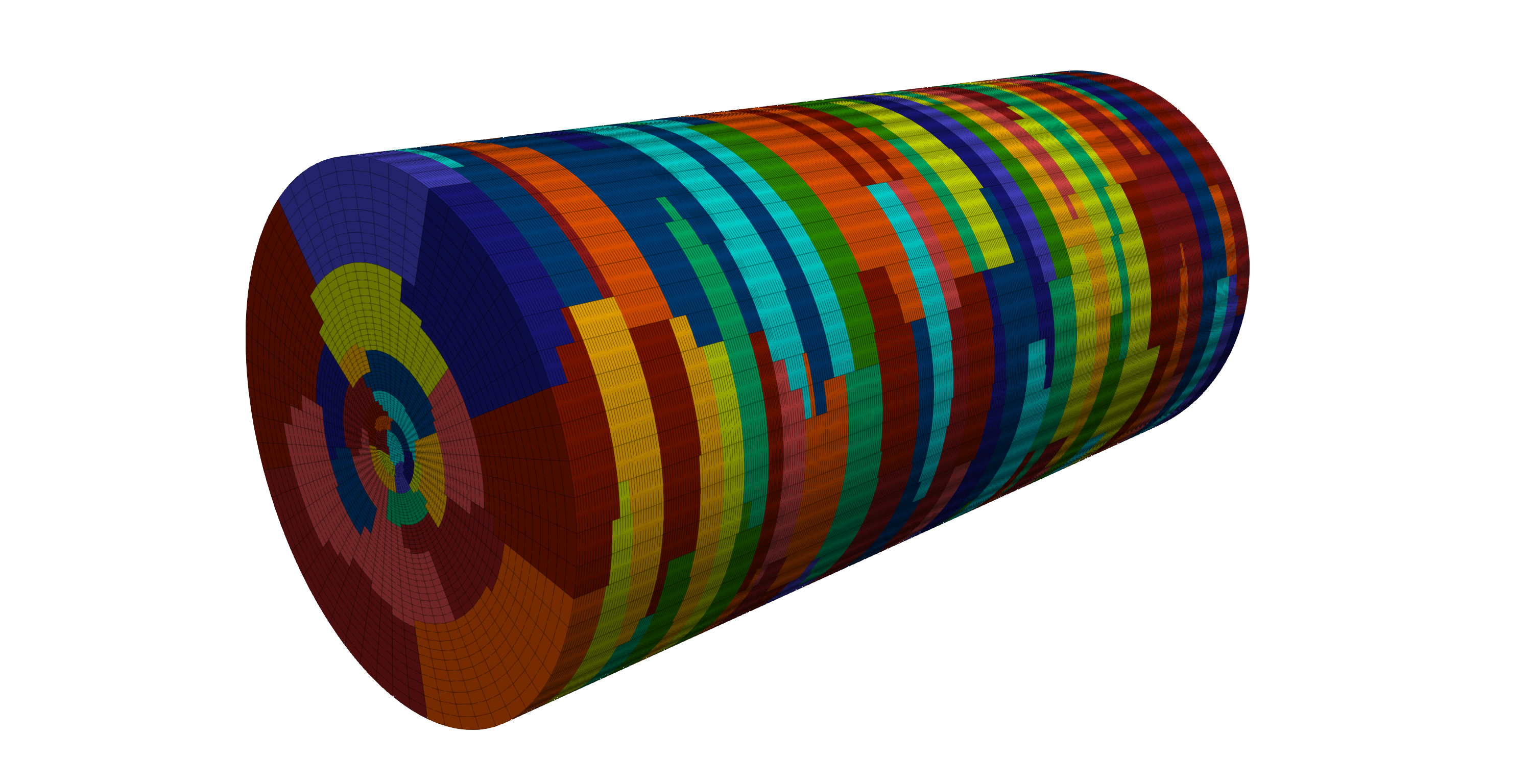}
\caption{Computational mesh for the flow around the unit sphere decomposed into 1024 subdomains.}
\label{fig:mesh}
\end{figure}

In this section, we present numerical results for the simulation of incompressible flow around a unit sphere.
We focus on multiple variants of the three-level BDDC method for the Poisson problem for pressure corrector~\eqref{eq:projection2} with a fixed matrix and time-dependent right-hand side vector.
The computational mesh of Taylor--Hood hexahedral elements leads to 1.4M unknowns for pressure, and it is decomposed into 1024 subdomains; see Fig.~\ref{fig:mesh}.
The time step is constant and set to 0.05 s.

The computations were performed on the \emph{Karolina} supercomputer at the IT4In\-no\-va\-tions National Supercomputing Centre in Ostrava, Czech Republic.
The computational nodes are equipped with two 64-core AMD 7H12 2.6 GHz processors, and 256 GB RAM.

\begin{figure}[!htp]
\centering
    \vskip -30mm
    \includegraphics[width=0.92\textwidth]{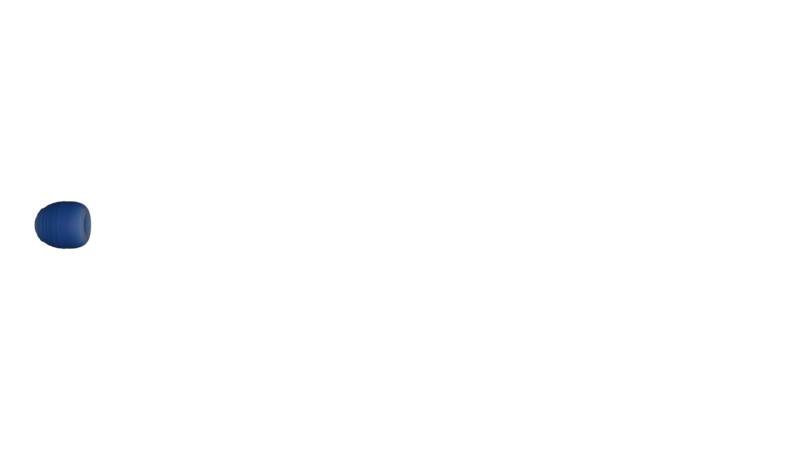} \\
    \vskip -60mm
    \includegraphics[width=0.92\textwidth]{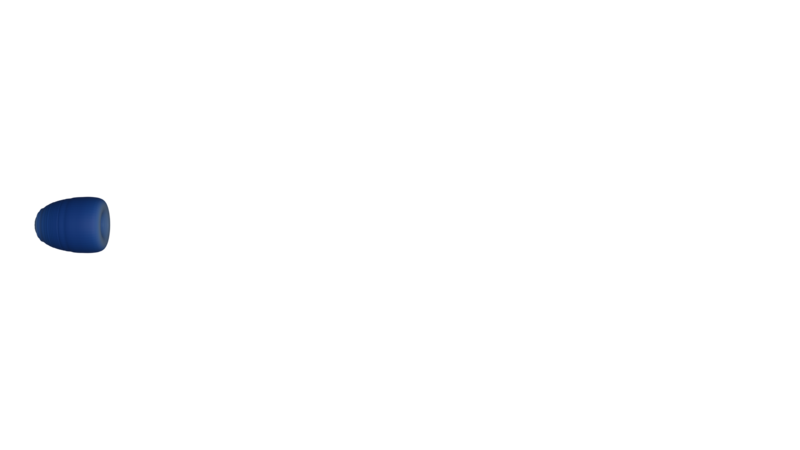} \\
    \vskip -60mm
    \includegraphics[width=0.92\textwidth]{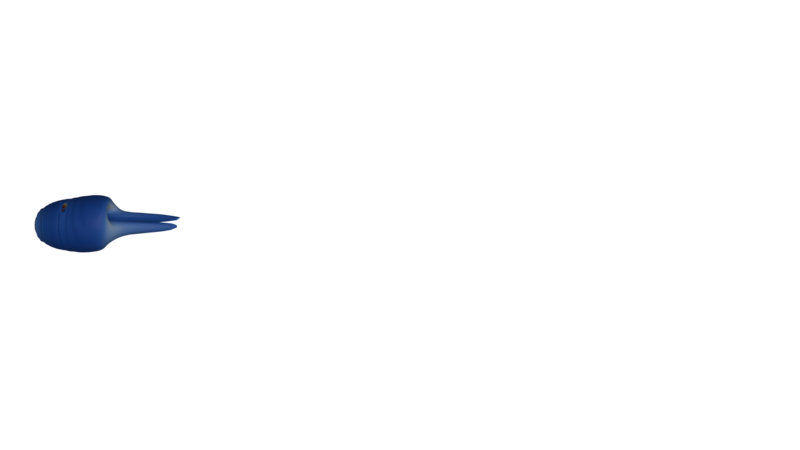} \\
    \vskip -60mm
    \includegraphics[width=0.92\textwidth]{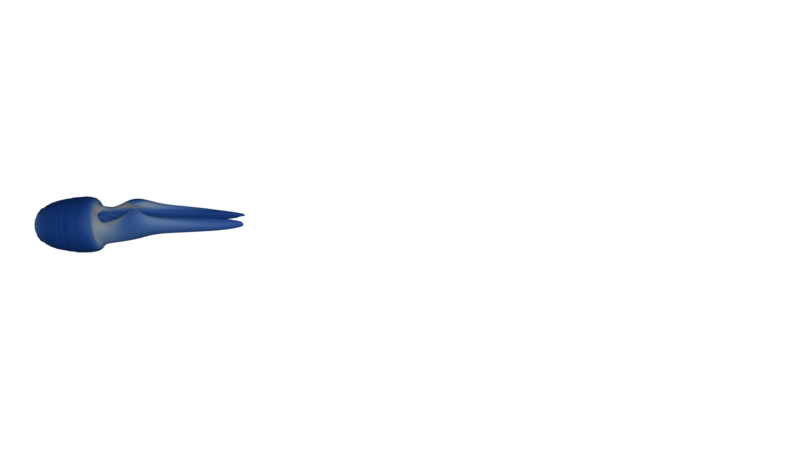} \\
    \vskip -60mm
    \includegraphics[width=0.92\textwidth]{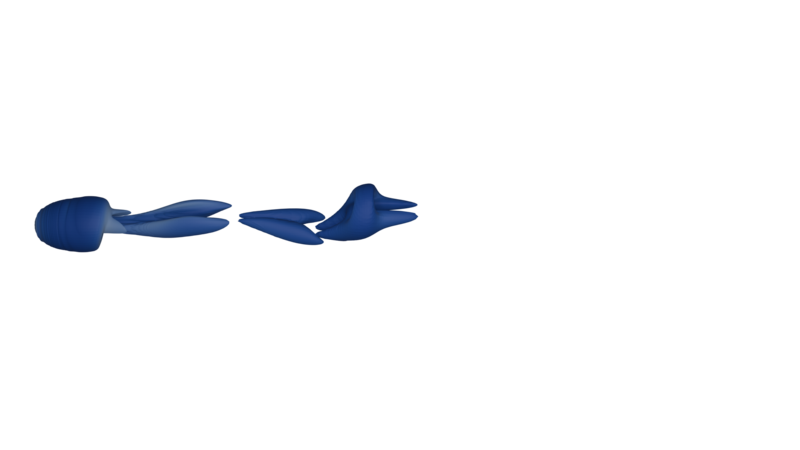} \\
    \vskip -60mm
    \includegraphics[width=0.92\textwidth]{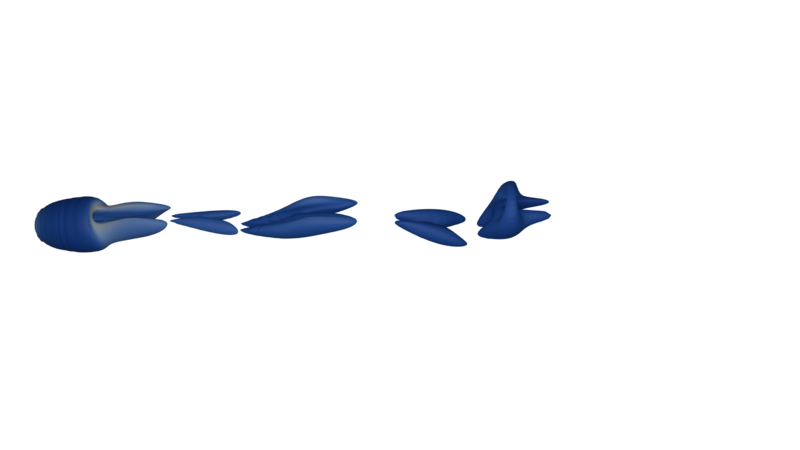} \\
    \vskip -60mm
    \includegraphics[width=0.92\textwidth]{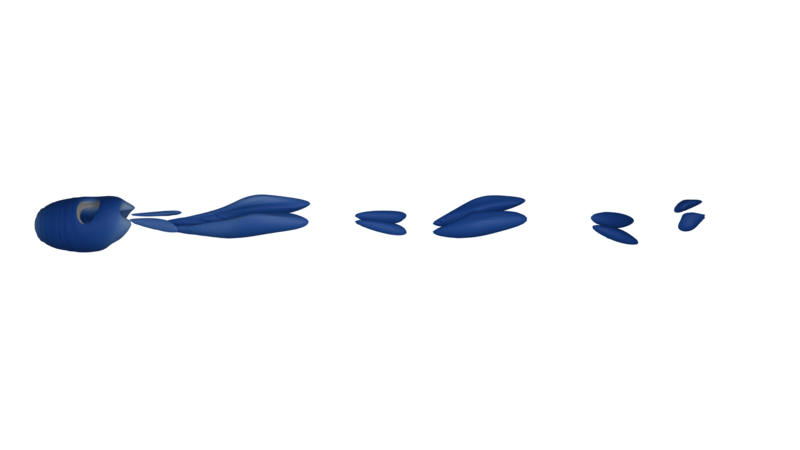} \\
    \vskip -60mm
    \includegraphics[width=0.92\textwidth]{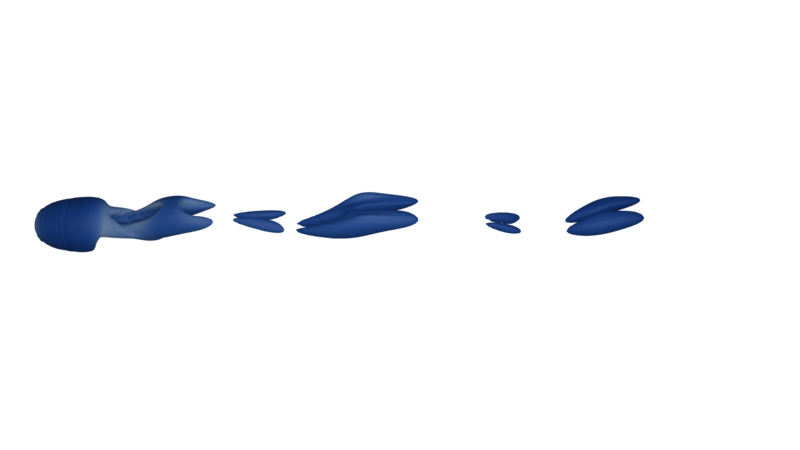} \\
    \vskip -60mm
    \includegraphics[width=0.92\textwidth]{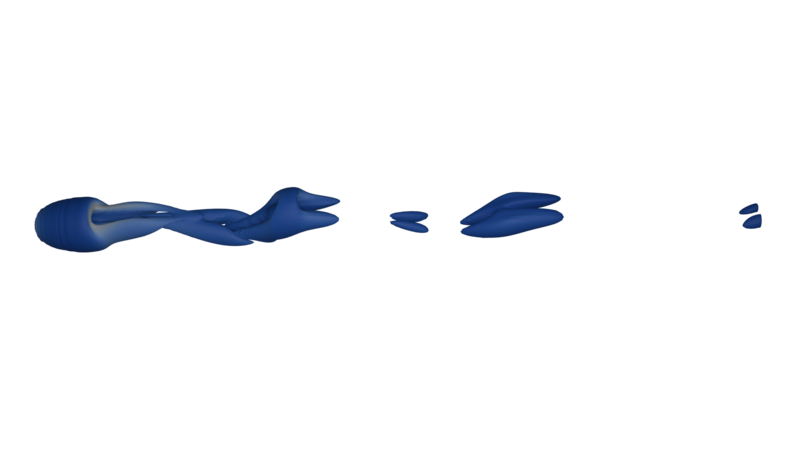} \\
    \vskip -60mm
    \includegraphics[width=0.92\textwidth]{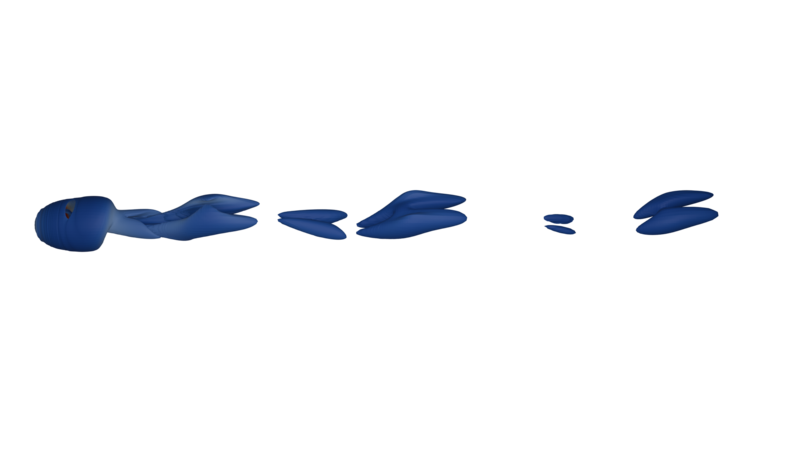} \\
    \vskip -60mm
    \includegraphics[width=0.92\textwidth]{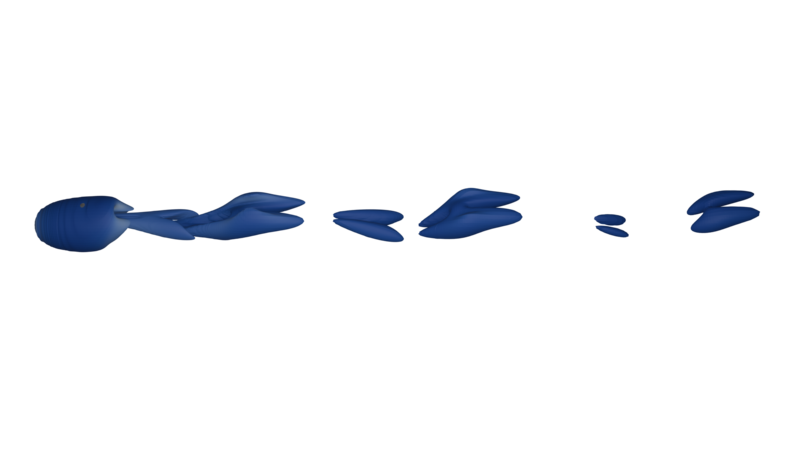} \\
    \vskip -60mm
    \includegraphics[width=0.92\textwidth]{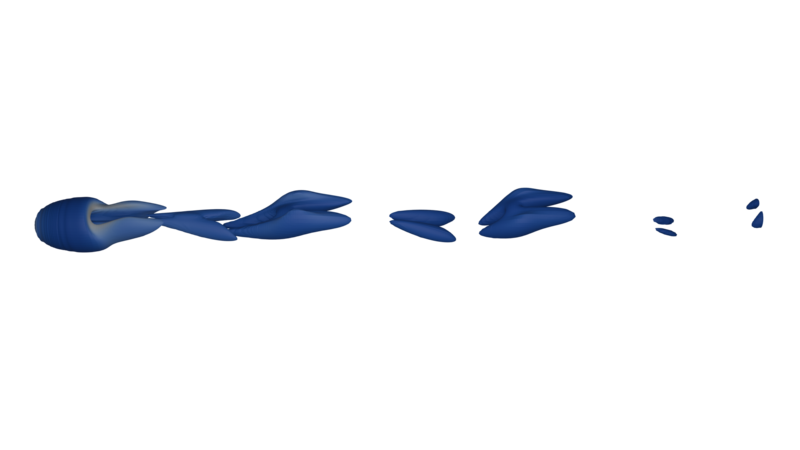} \\
    \vskip -60mm
    \includegraphics[width=0.92\textwidth]{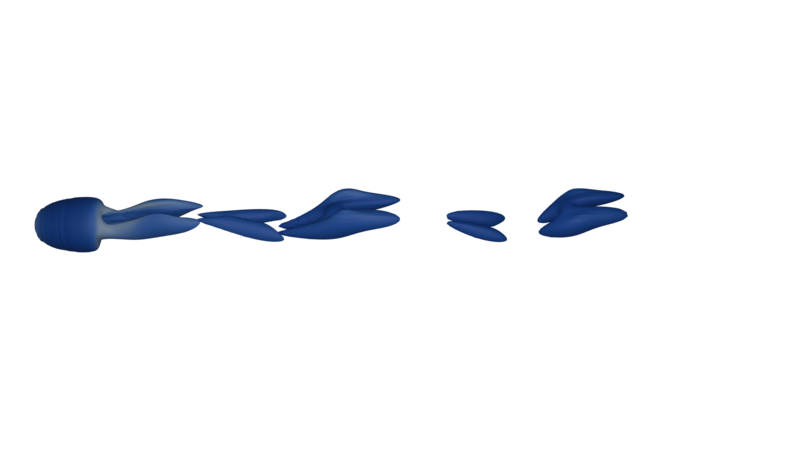} \\
    \vskip -60mm
    \includegraphics[width=0.92\textwidth]{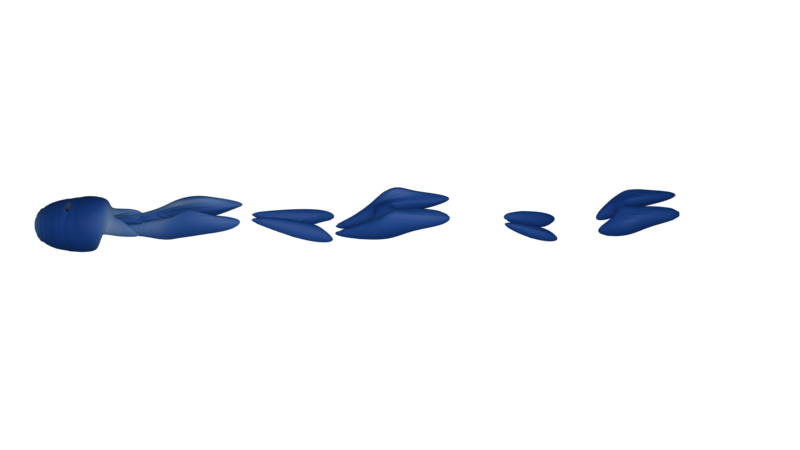}
    \vskip -30mm
\caption{Vortex structures behind a unit sphere at Re 300 at times 2.5, 25, 45, 50, 60, 65, 71, 75, 80, 100, 125, 150, 175, 200 s.}
\label{fig:vortices}
\end{figure}

\subsection{Effect of the stopping criterion}

\begin{figure}[htbp]
\centering
  \includegraphics[width=0.48\textwidth]{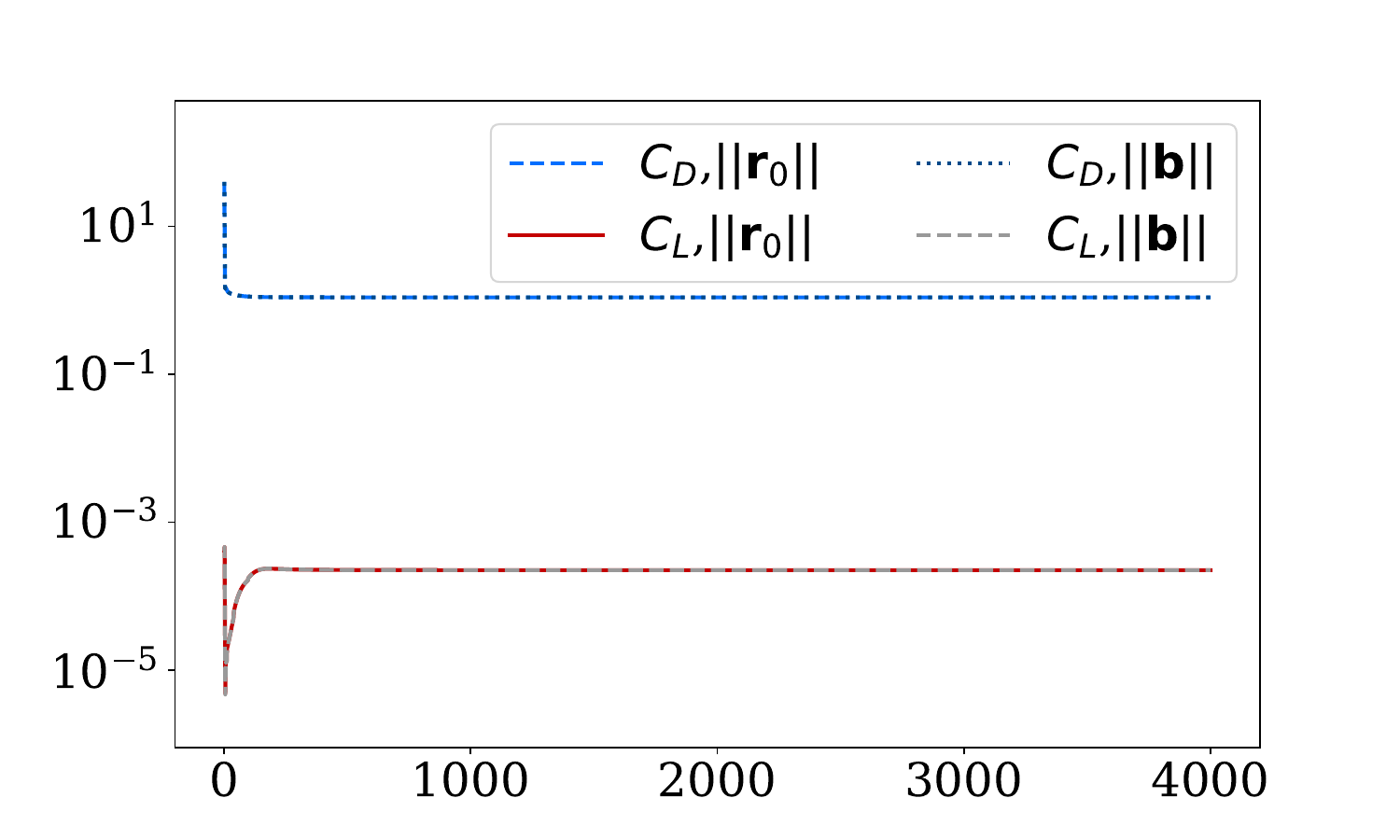}
  \includegraphics[width=0.48\textwidth]{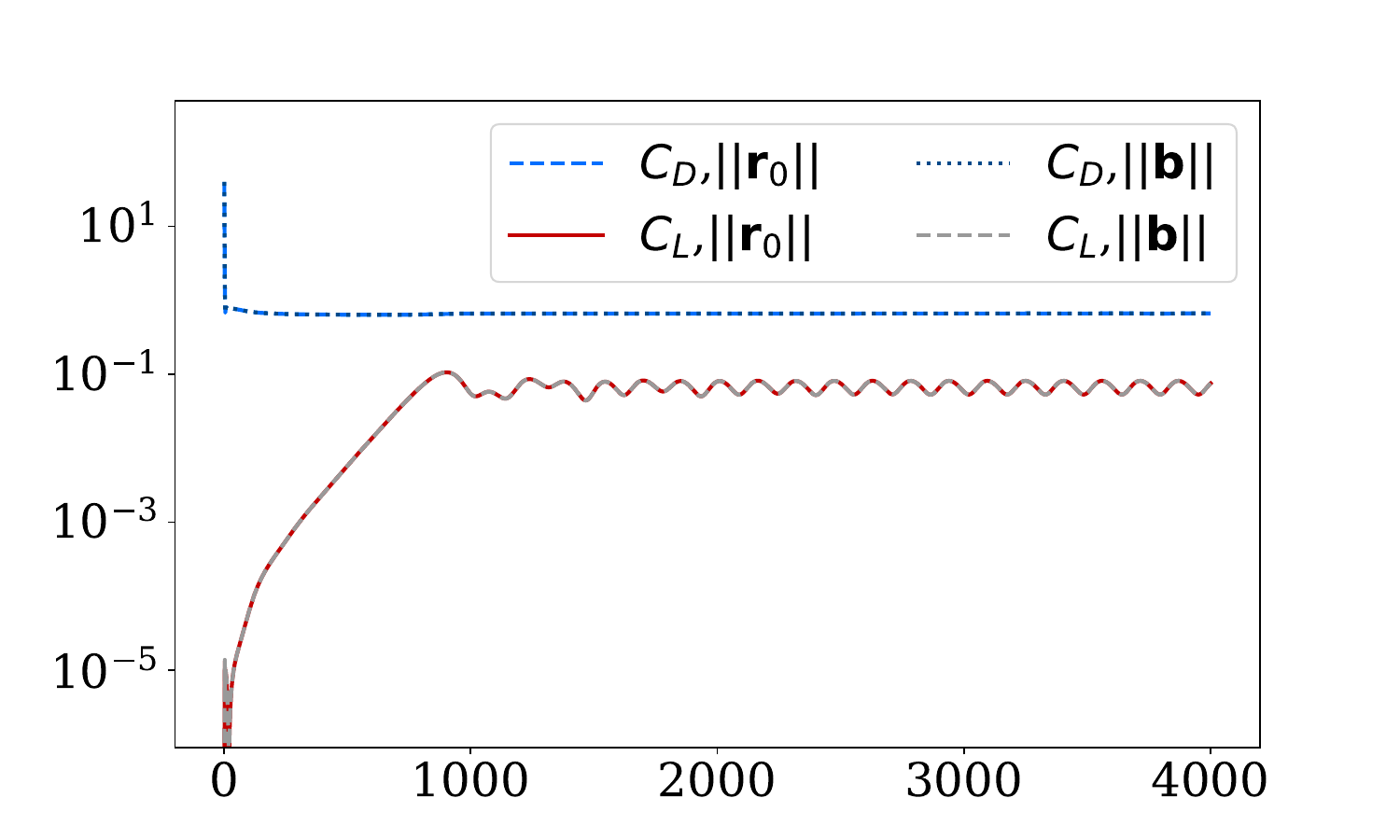}
\caption{Drag and lift coefficients over all time steps for Re = 100 (transient solution, left) and Re = 300 (periodic solution, right). 
The lines for different stopping criteria visually overlap, with relative difference in most of the iterations below $10^{-6}$.}
\label{fig:forces}
\end{figure}

Our main focus is the simulation of the flow at Reynolds number 300.
After an initial phase, the solution is periodic, see Fig.~\ref{fig:vortices} for several snapshots of vortex structures in different times and the right part of Fig.~\ref{fig:forces} for the evolution of the drag and lift coefficients in time.
In this experiment, we first focus on the effect of the choice of the initial guess and the stopping criteria,
and for that purpose, it is customary to perform also simulations for the case of Reynolds number 100.
For this lower Re, the solution is transient and converges to a steady state as can be deduced from the plot of the aerodynamic coefficients in Fig.~\ref{fig:forces} (left).
The vortex structure resembles the solution at the second snapshot of Fig.~\ref{fig:vortices}.
From the linear algebraic viewpoint, the right-hand sides in subsequent linear systems (\ref{eq:systems}) are becoming closer to each other, and the solution from the preceding time step is a progressively better approximation to the solution from the actual time step; recall that the system matrix and preconditioner remain the same in all time steps.
We are interested in the capability of the solver to exploit this fact and converge in a lower number of iterations.
The multilevel BDDC preconditioner uses the inverse cardinality (\emph{card}) as the interface weights.

\begin{figure}[htbp]
\centering
  \includegraphics[width=0.48\textwidth]{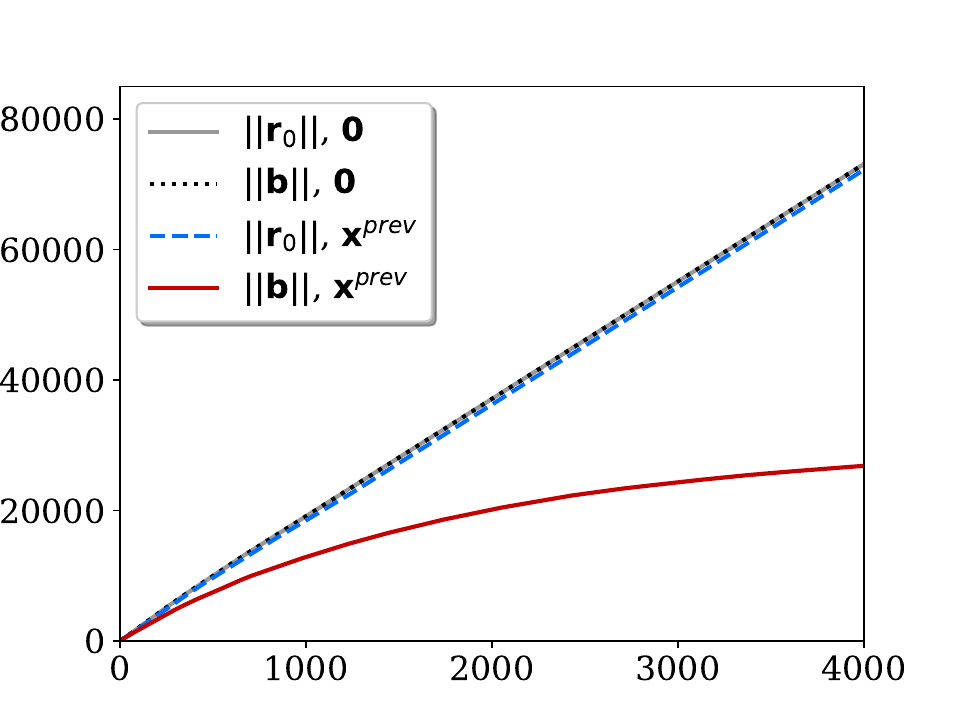}
  \includegraphics[width=0.48\textwidth]{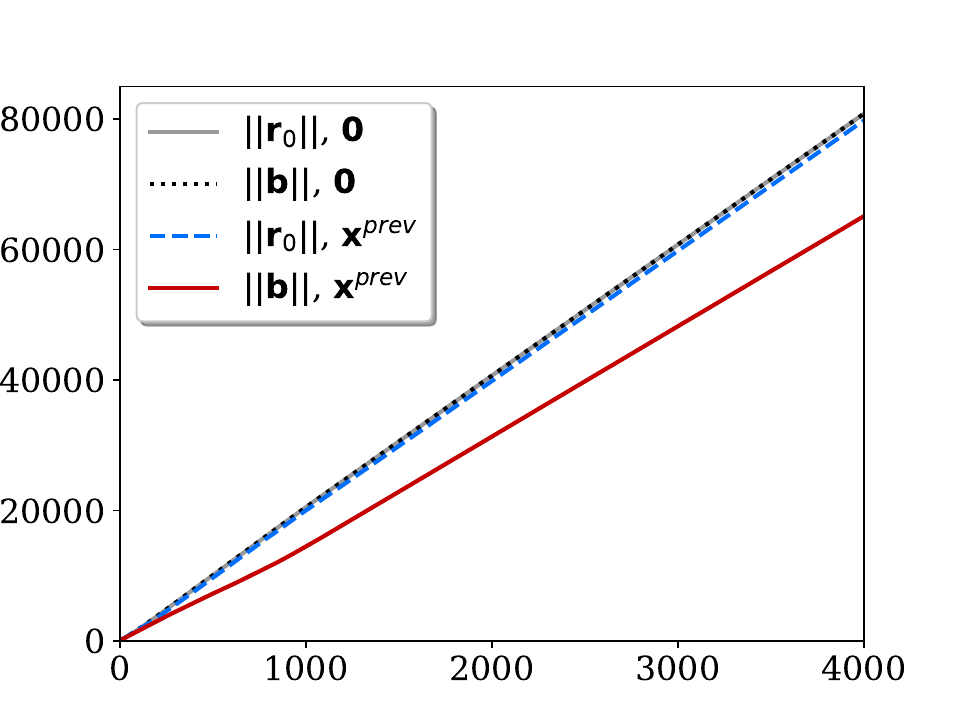}
\caption{Cumulative number of PCG iterations over all time steps for Re = 100 (transient solution, left) and Re = 300 (periodic solution, right).
Iterations terminated using ${\|\vr_{k}\|}/{\|\vr_{0}\|}$ (denoted by `$\|\vr_{0}\|$') and ${\|\vr_{k}\|}/{\|\vb\|}$ (`$\|\vb\|$') stopping criteria, initial guess taken as a zero vector (`$\vzero$') and as the approximate solution from the previous time step (`$\vx^\mathrm{prev}$').}
\label{fig:cumulative_iterations}
\end{figure}

The cumulative number of PCG iterations in 4000 time steps for the different configurations is presented in Fig.~\ref{fig:cumulative_iterations}, left.
We consider the choice of initial guess as the vector of zeros as well as the final approximation computed in the previous time step.
Two stopping criteria, \eqref{eq:stopping_residual} and \eqref{eq:stopping_rhs} are compared, both with $tol = 10^{-6}$.
The cumulative number illustrates the overall demand for the number of iterations of each approach over the course of the simulation.

From Fig.~\ref{fig:cumulative_iterations}, we can see that the lines corresponding to starting the solution always with the zero guess overlap
(when starting with the zero guess, the stopping criteria \eqref{eq:stopping_residual} and \eqref{eq:stopping_rhs} coincide as $\vr_0 = \vb$) and show a constant increase,
indicating that the number of iterations remains almost constant in each step.
When the solver
starts with the solution from the previous time step and the stopping criterion is based on the norm of the initial residual,
there is a marginal decrease in the number of iterations.
The desired behavior is obtained only for the stopping criterion~(\ref{eq:stopping_rhs}) based on the norm of the right-hand side vector, and the number of iterations decreases during the course of the time steps.

As illustrated in the left part of Fig.~\ref{fig:cumulative_iterations}, selecting the stopping criterion~(\ref{eq:stopping_rhs}) can significantly reduce the cumulative number of iterations for a simulation of the transient flow with Reynolds number 100,
in our case cutting the computational time by more than 50 percent.

\begin{figure}[htbp]
\centering
  \includegraphics[width=0.48\textwidth]{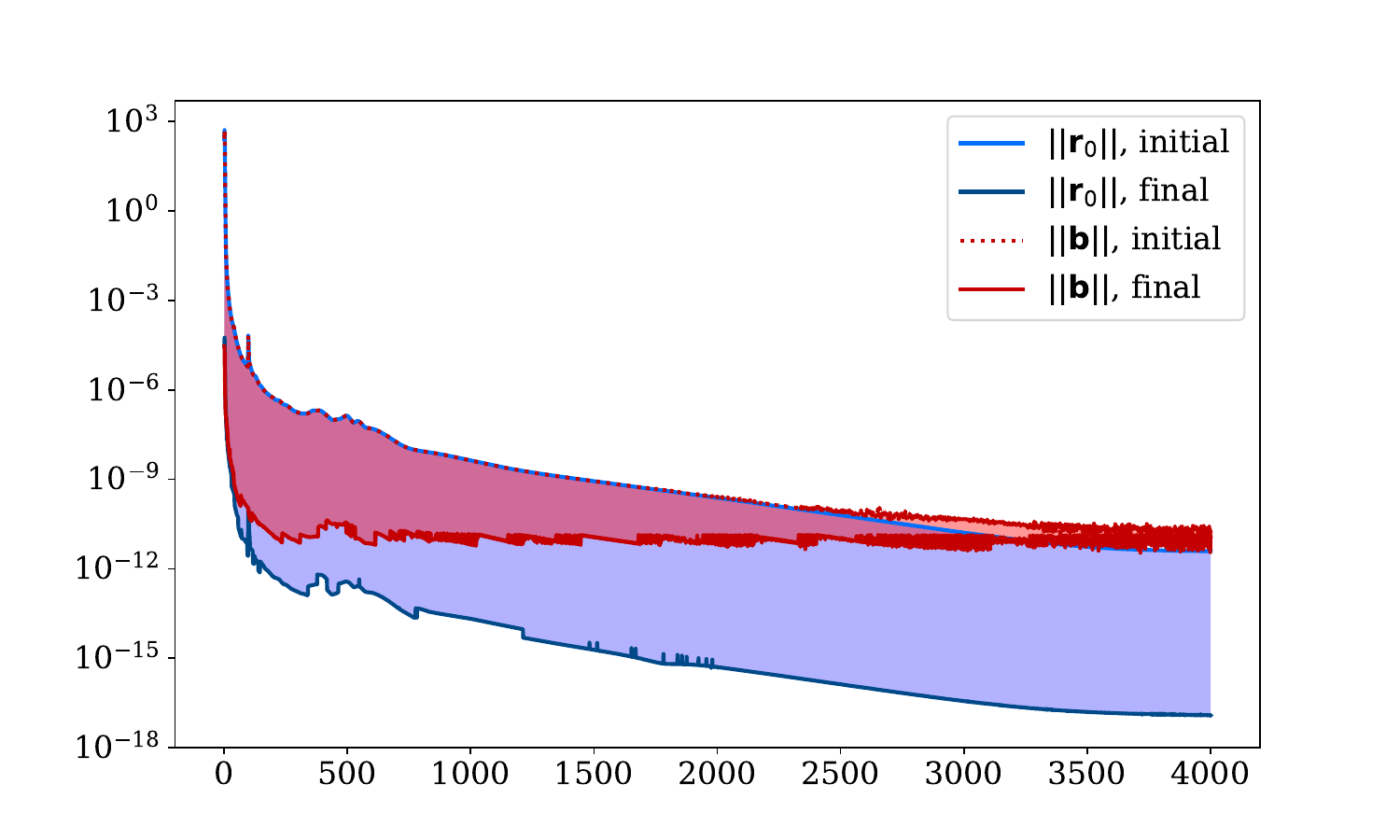}
  \includegraphics[width=0.48\textwidth]{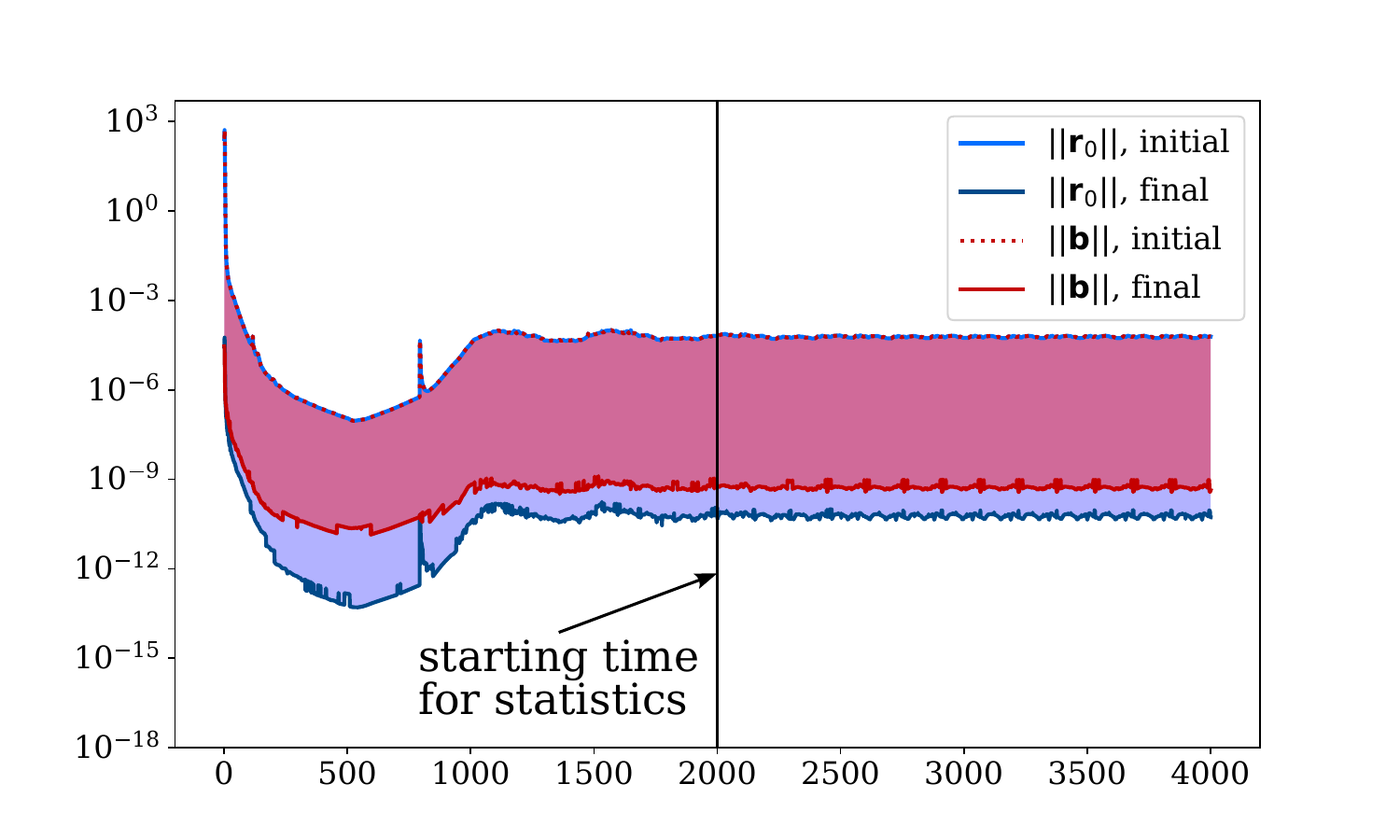}
\caption{Norm of the initial and final residual over all time steps for Re = 100 (transient solution, left) and Re = 300 (periodic solution, right).
Iterations terminated using ${\|\vr_{k}\|}/{\|\vr_{0}\|}$ (denoted as `$\|\vr_{0}\|$') and ${\|\vr_{k}\|}/{\|\vb\|}$ (`$\|\vb\|$') stopping criteria.
The time from which the statistics are computed is marked by a vertical line for Re = 300.}
\label{fig:residuals}
\end{figure}

In order to investigate this simulation closer, we present a plot of the residual norms at the beginning and at the end of each time step in Fig.~\ref{fig:residuals}. 
This clearly illustrates that while the residual norm of the initial guess computed from the previous approximation is decreasing in later time steps when the two consecutive solutions get closer to each other,
the stopping criterion based on $\|\vr_{0}\|$ leads to significant over-solving.
In most of the time steps, the residual norm of initial guesses nearly coincide showing that a softer stopping criterion is sufficient.
This is also confirmed when one compares the computed drag and lift coefficients in Fig.~\ref{fig:forces}.
The relative difference between the coefficients stays, apart from a few initial time steps, below~$10^{-6}$ indicating that the criterion~\eqref{eq:stopping_rhs} is indeed sufficient.

In our further experiments, we focus on the case with Re = 300, where a notable change of the regime results in a periodic behavior.
In this case, as shown in the right side of Fig.~\ref{fig:cumulative_iterations},
the choice of the initial guess and an appropriate stopping criterion alone does not lead to as substantial savings as for Re = 100, although the improvement is still considerable.
In the right part of Fig.~\ref{fig:residuals}, we can see that in the periodic regime the residual norm of the initial guess computed from the previous approximation does not decrease.
After an initial phase (approximately 1000 time steps) where the periodic solution is developed, the initial residual stagnates on a certain level (here around $6\times 10^{-5}$).
The change in the stopping criterion then allows for a larger residual norm of the computed approximation.
This softer criterion leads to a negligible difference in the resulting drag and lift coefficients; see the right panel in Fig.~\ref{fig:forces}.

Let us provide some more quantitative results for Re = 300.
In Table~\ref{tab:no_recycling_sphere3},
we report the minimum, maximum, and mean number of PCG iterations, the mean times for the whole linear solve and for one iteration,
and separately the number of iterations and time for the first step, which is excluded from the other statistics, as it includes the setup of the BDDC preconditioner.
For this table (and all tables hereafter), the statistics is computed from the time steps 2001 to 2100, i.e., after the periodic solution is developed, see Fig.~\ref{fig:residuals}.

\begin{table}[h]
\begin{center}
\footnotesize{
\begin{tabular}{c|c|cccc}
 stop. & init. & \#its.  & time [s] & step 1  \\
 crit. & guess &  min--max(avg) & step(1 iter.) & \#its. / time [s] \\ 
\hline\hline
\multirow{2}{*}{${\|\vr_{k}\|}/{\|\vr_{0}\|}$} & $\vzero$      & 20-20(20)           & 0.27(0.014)         & 20            /    0.79  \\ 
                                              & $\vx^\mathrm{prev}$ & 19-20(20.0)          & 0.28(0.014)          & 20            / 0.79 \\ 
\cline{1-2}
\multirow{2}{*}{${\|\vr_{k}\|}/{\|\vb\|}$}      & $\vzero$      & 20-20(20)          & 0.27(0.014)          & 20            / 0.76 \\
                                              & $\vx^\mathrm{prev}$ & 17-17(17)           & 0.24(0.014)         & 20            /    0.77  \\  
\hline
\end{tabular}
}
\caption{Varying stopping criterion and initial guess, Re = 300: Minimum (`min'), maximum (`max'), and mean (`avg') number of PCG iterations per time step, mean time for one step with mean time for one iteration in parentheses.
For the first step (system), the number of PCG iterations and the overall time including the construction of the BDDC preconditioner are reported.
\label{tab:no_recycling_sphere3}
}
\end{center} 
\end{table}

Table~\ref{tab:no_recycling_sphere3} confirms that for these simulations,
the most beneficial configuration is the variant with the stopping criterion~\eqref{eq:stopping_rhs} based on the norm of the right-hand side vector combined with the computed approximation from the previous time step taken as the initial guess for the new time step.
We therefore consider this choice of the initial guess and the stopping criterion in all experiments that follow.
As expected, the mean time for one PCG iteration and the number of iterations for the first step are consistent across all four cases. 

\subsection{Effect of Krylov subspace recycling}

In the next experiment, we focus on the strategies for recycling the Krylov subspace.
In particular, we consider four variants of the construction of the deflation basis from Section~{\ref{sec:recycling}} summarized in Table~\ref{tab:bases}.
For these simulations, the maximal size of the deflation basis is fixed at $R = 50$.
In addition to the results presented in Table~\ref{tab:no_recycling_sphere3}, we report the time step at which the norm of the Ritz values difference satisfied criterion~(\ref{eq:harmritzstopcrit}). After this step, the deflation basis was not updated.
Recall that the statistics is computed from the time steps 2001 to 2100.

\begin{table}[h]
\begin{center}
\begin{minipage}{0.9\textwidth}
\begin{description}
    \item[\textit{B1}] (first $R$ search vectors): The deflation basis is kept the same for all systems. It consists of the first $R$ search vectors for the first systems.
    \item[\textit{B2}] (last $R$ search vectors): The deflation basis consists of the last $R$ search vectors from solving the previous systems.
    \item[\textit{B3}] (eliminating $R$ smallest Ritz values): The basis contains approximations to~$R$ eigenvectors associated with the \emph{smallest} eigenvalues of the preconditioned matrix $\mM^{-1}\mA$.
                       These are computed using the harmonic Ritz approximation on the subspace generated by the deflation basis~(\ref{eq:harmonicRitz}) and the search vectors for the previous system.
                       It is updated in each time step until the criterion~(\ref{eq:harmritzstopcrit}) is satisfied.
    \item[\textit{B4}] (eliminating $R$ largest Ritz values): Similar to \textit{B3}, but the basis contains approximations associated with the $R$ \emph{largest} eigenvalues of the preconditioned matrix.
\end{description}
\end{minipage}
\end{center} 
\caption{Summary of the recycling strategies.}
\label{tab:bases}
\end{table}

\begin{table}[h]
\begin{center}
\footnotesize{
\begin{tabular}{c|ccccc}
 Recycl. & \#its.  & time [s] & step 1 & Ritz  \\
   strategy &  min--max(avg) & step(1 iter.) & \#its. / time [s] & conv. \\ 
\hline
\hline
{\textit{B1}} &  14-26(15.1) & 0.24(0.016)        & 20          / 0.82        & n/a\\  

{\textit{B2}} &  12-16(13.6) & 0.23(0.016)        & 20          / 0.84        & n/a\\  

{\textit{B3}} &  13-16(14.4) & 0.24(0.016)        & 20          / 1.76        & $>$100\\  

{\textit{B4}} & 13-16(13.4)  & 0.22(0.016)        & 20          / 1.70        & 46\\  
\hline

\end{tabular}
}
\caption{Results for different recycling strategies, Re = 300: Minimum (`min'), maximum (`max'), and mean (`avg') number of PCG iterations per time step, mean time for one time step with mean time for one \hl{PCG}~iteration in parentheses.
For the first step, the number of PCG iterations and the overall time including the construction of the BDDC preconditioner are reported.
The last column (`Ritz conv.') gives the index of the step when the criterion~\eqref{eq:harmritzstopcrit} is met, i.e., when the Ritz values converged.
The size of the deflation basis is set to $R=50$.
\label{tab:recycling_cases_sphere3}
}
\end{center} 
\end{table}

From Table~\ref{tab:recycling_cases_sphere3} we observe that all strategies reduce the number of PCG iterations, which was originally 17.
The time of single PCG iteration is the same, up to two valid digits, for all the recycling strategies,
and it increased by~14~\% in comparison with the variant without deflation.
Note that \textit{B1} and \textit{B2} do not require recomputing the deflation basis using~\eqref{eq:harmonicRitz},
but the cost of recomputing the basis is negligible with respect to other operations.
Overall, the savings in computational time per step are not very satisfactory for  \textit{B1} and \textit{B3}.
Due to the slow convergence of the smallest Ritz values, the deflation basis is recomputed in each time step when using~\textit{B3}.

The best results were obtained for the strategy \textit{B4}, which uses the $R$~vectors corresponding to the largest Ritz values. This choice requires, on average, the lowest number of PCG iterations, and the lowest computational time for one time step.
In this setting, also the strategy \textit{B2}, which is a bit simpler to implement, can also be recommended.
Following this experiment, we restrict ourselves to the recycling strategy \textit{B4} for further experiments.

In the next experiment, we vary the size of the deflation basis~$R$. The results are presented in Table~\ref{tab:recycling_base_sphere3}.
As expected, increasing the size of the deflation basis~$R$ reduces the mean number of iterations, but increases the time required for each of them. 
Choosing a proper size of the deflation basis can be a difficult task and, in practice, this might be done based on memory limitations, i.e., choosing the size as the maximal number of the vectors that can be stored. This experiment reveals, however, that similar timings can be obtained for a range of sizes.
The choice $R=50$, used in the previous experiment, is one with nearly the best computational time for a single time step.
To provide a fair comparison, we keep the size of the deflation basis $R=50$ also for the following numerical experiments.

\begin{table}[h]
\begin{center}
\footnotesize{
\begin{tabular}{c|ccccc}
 basis size  & \#its.  & time [s] & step 1 & Ritz  \\
   $R$   &  min--max(avg) & step(1 iter.) & \#its. / time [s] & conv. \\ 
\hline
\hline
  25             & 16-17(16.2)           & 0.24(0.015)          & 20           / 1.88          & 12\\

  30             & 13-16(13.9)           & 0.22(0.016)          & 20            / 1.75          & 15\\
  
  35             & 13-16(13.2)           & 0.21(0.016)          & 20            / 1.81          & 17\\

  50             & 13-16(13.4)           & 0.22(0.016)          & 20            / 1.70          & 46\\

  100            & 9-16(11.7)            & 0.24(0.020)          & 20            / 1.78          & $>$100\\

  200            & 8-16(10.2)             & 0.23(0.022)          & 20            / 1.97          & $>$100\\

  400            & 8-16(10.3)             & 0.23(0.022)          & 20            / 2.07          & $>$100\\
\hline
\end{tabular}
}
\caption{Varying the maximal size of the recycling basis $R$ in recycling strategy \textit{B4}:
Minimum (`min'), maximum (`max'), and mean (`avg') number of PCG iterations per time step, mean time for one time step with mean time for one \hl{PCG}~iteration in parentheses.
For the first step, the number of PCG iterations and the overall time including the construction of the BDDC preconditioner are reported.
The last column (`Ritz conv.') gives the index of the step when the criterion \eqref{eq:harmritzstopcrit} is met.
\label{tab:recycling_base_sphere3}
}
\end{center} 
\end{table}

\subsection{Effect of adaptive BDDC preconditioner}

In the next set of numerical experiments, we perform computations for variants of the adaptive selection of coarse degrees of freedom in the BDDC preconditioner on top of the recycling of the Krylov subspace with $R = 50$. For preliminary results without recycling and with recycling using strategy \textit{B1}, see our paper \cite{Hanek-2023-AMB}. In the current experiment, we test several values of the prescribed target value $\tau$ described in Section~\ref{sec:adaptiveBDDC}. For a smaller $\tau$, more eigenvectors are used in the construction of the coarse problem. This reduces the number of iterations, but again, each iteration gets more expensive due to a larger coarse problem. The method is combined with two types of weights in BDDC (i.e., matrices $\mD_i$ introduced in Section~\ref{sec:multilevelBDDC}), namely the scaling based on the cardinality (\emph{card}) and based on the diagonal stiffness (\emph{diag}).

\begin{table}[h]
\begin{center}
{\footnotesize
\begin{tabular}{c|c|ccccccc}
 \multicolumn{2}{c|}{adaptivity} & \#its.  & time [s] & step 1 & Ritz  \\
   weights  & $\tau$  &  min--max(avg) & step(1 iter.) & \#its. / time [s] & conv. \\ 
\hline\hline
 \multirow{4}{*}{card}                    & 3.5    & 9-12(11.3)            & 0.20(0.017)          & 14            / 17.85         & 73\\
                                          & 3.0    & 8-12(10.3)             & 0.19(0.018)          & 14            / 18.23         & 74\\
                                          & 2.5    & 8-11(10.0)             & 0.19(0.019)          & 14            / 18.85         & 69\\
                                          & 2.0    & 8-11(8.7)             & 0.25(0.025)          & 14            / 21.69         & 76\\ \hline
                    \multirow{4}{*}{diag} & 3.5    & 10-12(11.0)             & 0.19(0.017)          & 15            / 18.14         & 33\\
                                          & 3.0    & 10-12(10.2)             & 0.18(0.018)          & 15            / 18.41         & 39\\
                                          & 2.5    & 10-12(10.1)             & 0.19(0.019)          & 15            / 19.78         & 33\\
                                          & 2.0    & 9-12(10.0)             & 0.20(0.020)          & 16            / 19.80            & 30\\
\hline
\end{tabular}
}
\caption{Varying the adaptive coarse space in BDDC by changing the threshold on eigenvalues for selecting eigenvectors for the coarse problem $\tau$: Minimum (`min'), maximum (`max'), and mean (`avg') number of PCG iterations per time step, mean time for one time step with mean time for one \hl{PCG}~iteration in parentheses.
The size of the deflation basis $R$ = 50.
For the first step, the number of PCG iterations and the overall time including the construction of the BDDC preconditioner are reported.
The last column (`Ritz conv.') gives the index of the step when the criterion \eqref{eq:harmritzstopcrit} is met.
\label{tab:recycling_case2_sphere3}
}
\end{center} 
\end{table}

From the results in Table~\ref{tab:recycling_case2_sphere3}, we can conclude that the optimal threshold for eigenvalues $\tau$ is 3.0 for both types of weights, in terms of the mean time for the linear solver.
The results indicate that the weights based on diagonal stiffness (\emph{diag}) perform slightly better than those based on cardinality (\emph{card}).
Consequently, \emph{diag} weights are used in further experiments of this paper.

It is important to note that in our current implementation, the adaptive BDDC method introduces a substantial additional cost to the preconditioner setup, as it requires solving local eigenvalue problems~(\ref{eq:eig-matrix}). By comparing the average time per one time step in Tables~\ref{tab:recycling_cases_sphere3} and \ref{tab:recycling_case2_sphere3}, we observe that approximately 350 time steps are needed to amortize this initial cost.

\begin{figure}[htbp]
\centering
   \ding{56} recycling \quad \ding{56} adaptivity \\
   \includegraphics[width=0.85\textwidth]{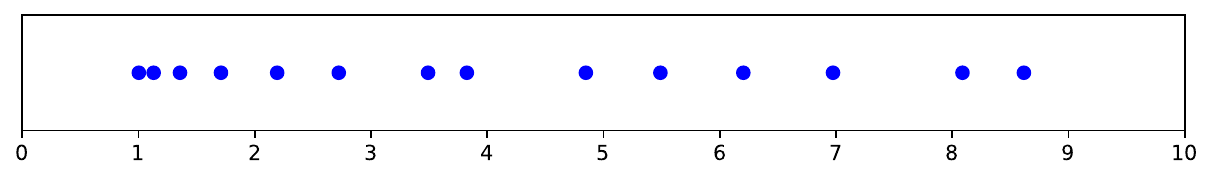} \\
   \ding{52} recycling \quad \ding{56} adaptivity \\
   \includegraphics[width=0.85\textwidth]{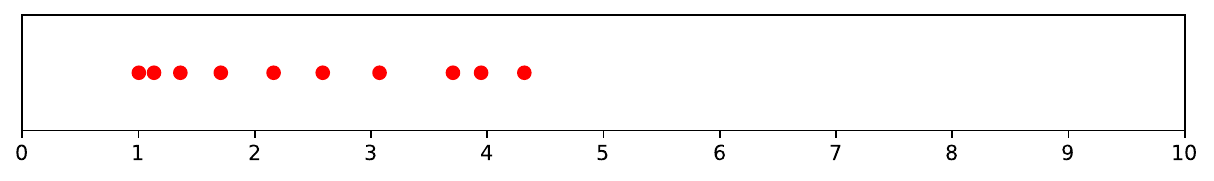}\\
   \ding{56} recycling \quad \ding{52} adaptivity \\
   \includegraphics[width=0.85\textwidth]{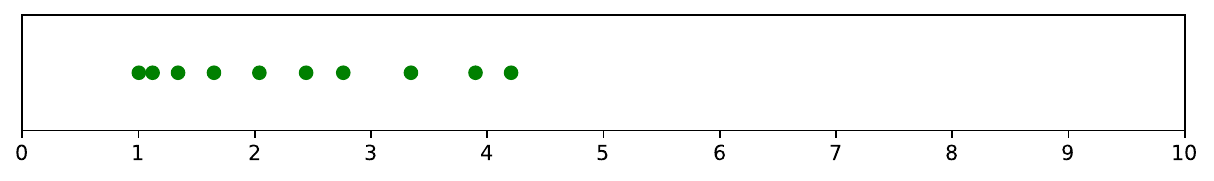}\\
   \ding{52} recycling \quad \ding{52} adaptivity \\
   \includegraphics[width=0.85\textwidth]{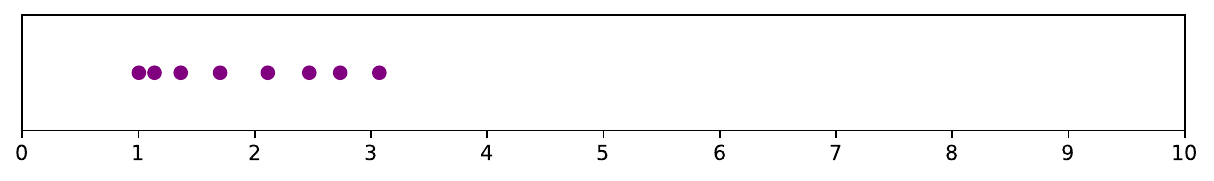}
\caption{Ritz values approximating the eigenvalues of the preconditioned system. From top to bottom: no recycling nor adaptive BDDC, Krylov subspace recycling with $R=50$, no recycling and adaptive BDDC with $\tau = 3.0$, and Krylov subspace recycling with $R=50$ and adaptive BDDC with $\tau = 3.0$.}
\label{fig:eigen_values}
\end{figure}

Table~\ref{tab:recycling_sphere3_sum} presents a summary of the results, comparing the best variant for each approach. The results show a 25\% speed-up in computational time for one time step for the variant that incorporates both recycling and adaptivity compared to the one without these approaches. To further understand this behavior, Figure~\ref{fig:eigen_values} displays the Ritz values approximating the spectrum of the preconditioned operator $\mM^{-1}\mA$, or $\mQ\mM^{-1}\mA$ when deflation is used. The figure illustrates how Krylov subspace recycling and the adaptive BDDC preconditioner push the upper part of the spectrum closer to one, thereby reducing the condition number and accelerating convergence.

\begin{table}[h]
\begin{center}
\footnotesize{
\begin{tabular}{cc|ccccc}
 \multicolumn{2}{c|}{variant} & \#its.  & time [s] & step 1  \\
  recycling & adaptivity &  min--max(avg) & step(1 iter.) & \#its. / time [s] \\ 
\hline\hline
  \ding{56} & \ding{56}  & 17-17(17)           & 0.24(0.014)         & 20            / \phantom{4}0.77 \\

  $R = 50$ & \ding{56}   & 13-16(13.4)            & 0.22(0.016)         & 20            / \phantom{4}1.70 \\

  \ding{56} & diag, $\tau = 3.0$   & 12-13(12.8)           & 0.19(0.015)         & 15            / 18.79 \\

  $R = 50$ & diag, $\tau = 3.0$ & 10-12(10.2)             & 0.18(0.018)         & 15            / 18.41 \\
\hline

\end{tabular}
}
\caption{Summary of the best results from Tables~\ref{tab:no_recycling_sphere3}--\ref{tab:recycling_case2_sphere3} for different acceleration strategies.
Minimum (`min'), maximum (`max'), and mean (`avg') number of PCG iterations per time step, mean time for one time step with mean time for one \hl{PCG}~iteration in parentheses, and the number of PCG iterations (`step 1 \#its.'), and time for the first time step. We consider the initial guess given by the computed solution to the previous system and the stopping criterion~\eqref{eq:stopping_rhs}.
}
\label{tab:recycling_sphere3_sum}
\end{center} 
\end{table}

\medskip

To investigate the efficiency of the acceleration strategies studied for the smaller mesh size, we finally consider a larger computational mesh, consisting of 15.6 million unknowns for pressure solved on 4096 subdomains and processor cores.
We monitor the same metrics of iterations and times as for the smaller problem.
Benchmark simulations are carried out again in four different modes:
without Krylov subspace recycling, with recycling, with adaptive BDDC, and with a combination of both Krylov subspace recycling and adaptive BDDC.
The results are summarized in Table~\ref{tab:recycling_with_adaptivity_sphere4}.

\begin{table}[h]
\begin{center}
\footnotesize{
\begin{tabular}{cc|ccccc}
 \multicolumn{2}{c|}{variant} & \#its.  & time [s] & step 1  \\
  recycling & adaptivity &  min--max(avg) & step(1 iter.) & \#its. / time [s] \\ 
\hline\hline
  \ding{56} & \ding{56}  & 37-38(37.7)           & 3.94(0.104)         & 44            / \phantom{4}7.00 \\

  $R = 50$ & \ding{56}   & 26-30(26.1)            & 2.93(0.113)         & 44           / \phantom{4}8.17 \\

  \ding{56} & diag, $\tau = 3.0$   & 20-22(20.6)           & 2.53(0.122)         & 24            / 114.06 \\

  $R = 50$ & diag, $\tau = 3.0$ & 16-19(16.8)             & 2.21(0.132)         & 24            / 113.62 \\
\hline

\end{tabular}
}
\caption{Larger computational mesh with 15.6 million pressure unknowns and 4096 subdomains: analogy to the results from Table~\ref{tab:recycling_sphere3_sum}.}
\label{tab:recycling_with_adaptivity_sphere4}
\end{center} 
\end{table}

The trends observed for the smaller problem are confirmed for the larger problem (compare Table~\ref{tab:recycling_sphere3_sum} and Table~\ref{tab:recycling_with_adaptivity_sphere4}).
When recycling, adaptivity, and their combination are introduced, the results align with those from the smaller problem, but with the differences being even more pronounced.
The combination of adaptive BDDC and Krylov subspace recycling proves to be the most effective in reducing the iteration count and, more importantly, the computational time for the linear solver.
In particular, these acceleration techniques lead to reducing the computational time per one time step by more than~40~\%.

When comparing the times reported in Tables~\ref{tab:recycling_sphere3_sum} and \ref{tab:recycling_with_adaptivity_sphere4}, we observe a noticeable increase for the larger problem.
This behavior is partly in line with the design of our experiments: since the second problem employs larger subdomains, weak scalability is not expected.
Furthermore, by constraining the size of the local problems for the velocity unknowns, the corresponding local problems for the pressure corrector become very small, and hence communication overhead dominates the computational time.

For completeness, we also performed experiments with time step sizes that were 10$\times$ larger and 10$\times$ smaller. This variation had a negligible effect on the convergence of the Poisson problem for the pressure corrector, and thus the conclusions of the previous experiments remain valid.

\section{\hl{Conclusions}}

We have focused on speeding up the solution of the Poisson problem for pressure corrector in time-depen\-dent incompressible flows problems. From the linear algebraic viewpoint, the problem leads to solving a sequence of algebraic systems with a constant symmetric positive-definite matrix and a right-hand side vector that changes in each time step.
A baseline for our study was provided by the three-level BDDC preconditioner within PCG. Three algorithmic components have been tested, and their effects on the time-to-solution have been evaluated.

We first studied the influence of the initial guess and the stopping criterion. To benefit from a good initial approximation, the stopping criterion must be independent of the initial guess. We therefore used a criterion normalized by the norm of the right-hand side. This allows the Krylov method to reduce the iteration count as the initial guess improves (e.g., when reusing the previous time-step solution), resulting in more than 50\% savings in computational time in the transient regime for Reynolds number Re = 100.

The second part of the study focused on deflation in PCG combined with Krylov subspace recycling. We tested four strategies of constructing the deflation basis: two based on storing previous search directions and two based on harmonic Ritz approximations of eigenvectors. Unlike the standard approach, which targets the smallest eigenvalues, we exploit the fact that preconditioners such as multilevel BDDC yield spectra bounded below by one and therefore we propose constructing the deflation basis from approximations of eigenvectors associated with the largest eigenvalues. This shifts down the upper end of the spectrum and, despite the higher per-iteration cost, reduces the total number of iterations and overall computational effort.

The final acceleration component was the adaptive selection of the BDDC coarse space. Adaptive BDDC tunes the preconditioner strength by solving local eigenvalue problems on subdomain interfaces and constructing a coarse space that reduces the condition number according to a threshold~$\tau$. Smaller values of~$\tau$ introduce more coarse degrees of freedom, decreasing the iteration count but increasing the per-iteration cost. Because the setup of the adaptive BDDC preconditioner is relatively expensive, this approach is particularly advantageous for longer simulations.

We demonstrated the synergy of three acceleration techniques for time-dependent problems. Using an appropriate stopping criterion, the previous time-step solution as the initial guess, deflation in PCG, and adaptive coarse space construction in the multilevel BDDC preconditioner led to substantial reductions in computational time. In the transient regime (Reynolds number Re = 100), more than 50\% of the time was saved by the stopping criterion alone. For the periodic regime with Re = 300, we achieved savings of about 25\% for the smaller problem with 1.4M pressure unknowns and over 40\% for the larger problem with 15.6M unknowns. These results indicate that Krylov subspace recycling and adaptive coarse spaces become increasingly beneficial for more challenging problems requiring many PCG iterations. While the individual components are well known, their synergistic combination has not been previously reported in the literature.

Finally, we note that both the size of the deflation basis in deflated PCG and the threshold~$\tau$ in adaptive BDDC admit optimal values. Please note that the goal of the paper was not to show the optimal setup of the parameters.
Instead, we identified them approximately and showed that they are not sharp, as comparable speed-ups are obtained over a range of parameter values. This indicates that the proposed approach is robust and broadly applicable to simulations of unsteady incompressible flows.

\section*{Acknowledgements}
This work was supported by the Czech Science Foundation through Grant No.~23-06159S.
It was also supported by the Institute of Mathematics of the Czech Academy of Sciences (RVO:67985840).
The computational time on the systems at IT4Innovations was provided thanks to the support by the Ministry of Education, Youth and Sports of the Czech Republic through the e-INFRA~CZ (ID:90254).
The authors would like to thank Žofia Machová for her help with rendering the vortex structures.
J.~Papež and J.~Šístek are members of the Nečas Center for Mathematical Modeling.

\section*{Data availability}
Data will be made available on request.

\section*{Declaration of competing interest}
The authors declare that they have no known competing financial interests or personal relationships that could have appeared to influence the work reported in this paper.


\end{document}